\def\a             {\alpha}
\def\Ad            {{\mathrm{Ad}}}
\newcommand\alt[3] {\tilde{\a}_{{#1};{#3}}^{{#2}}}

\def\be            {\begin{equation}}

\def\bbR           {\mathbb{R}}

\def\bbZ           {\mathbb{Z}}
\def\bfe           {{\bf1}}
\def\can           {\gamma}
\def\canr          {\theta}

\def\cA            {{\mathcal{A}}}

\def\cH            {{\mathcal{H}}}

\def\cK            {{\mathcal{K}}}

\def\cO            {{\mathcal{O}}}

\def\cZ            {{\mathcal{Z}}}
\newcommand\co[1]  {\overline{{#1}}}

\def\dim           {{\mathrm{dim}}}
\newcommand\del[2] {\delta_{{#1},{#2}}}
\newcommand\DS[1]  {\displaystyle{#1}}
\def\E             {{\mathrm{e}}}
\def\ee            {\end{equation}}
\def\End           {{\mathrm{End}}}
\def\eps           {\varepsilon}
\def\epsr          {\varepsilon_{\mathrm{r}}}
\newcommand\erf[1] {Eq.\ (\ref{#1})}

\def\ext           {{\mathrm{ext}}}

\def\Hom           {{\mathrm{Hom}}}
\def\I             {{\mathrm{i}}}
\def\id            {{\mathrm{id}}}

\def\la            {\lambda}
\def\lan           {\langle}

\def\MXM           {{}_M {\cal X}_M}
\def\MXMa          {{}_M^{} {\cal X}_M^\a}
\def\MXMo          {{}_M^{} {\cal X}_M^0}
\def\MXMop         {{}_{M_+}^{}\!\! {\cal X}_{M_+}^0}
\def\MXMom         {{}_{M_-}^{}\!\! {\cal X}_{M_-}^0}
\def\MXMopm        {{}_{M_\pm}^{}\!\! {\cal X}_{M_\pm}^0}
\def\MXMp          {{}_M^{} {\cal X}_M^+}
\def\MXMm          {{}_M^{} {\cal X}_M^-}
\def\MXMpm         {{}_M^{} {\cal X}_M^\pm}
\def\MXMpp         {{}_{M_+}^{} {\cal X}_{M_+}^+}
\def\MXMmm         {{}_{M_-}^{} {\cal X}_{M_-}^-}
\def\MXMpmpm       {{}_{M_\pm}^{} {\cal X}_{M_\pm}^\pm}
\def\NXN           {{}_N {\cal X}_N}
\def\NXNd          {{}_N^{} {\cal X}_N^{\mathrm{deg}}}

\def\om            {\omega}

\def\op            {{\mathrm{opp}}}

\def\ran           {\rangle}

\def\rmv           {{\mathrm{v}}}
\def\rms           {{\mathrm{s}}}
\def\rmc           {{\mathrm{c}}}

\def\sig           {\sigma}

\def\SLZ           {{\mathit{SL}}(2;\bbZ)}

\def\SOn           {{\mathit{SO}}(n)}
\def\SOs           {{\mathit{SO}}(16\ell)}

\def\SUd           {{\mathit{SU}}(3)}

\def\SUn           {{\mathit{SU}}(n)}

\def\SUz           {{\mathit{SU}}(2)}

\newcommand\tmat[1]{{}^{{\rm t}} {#1}}
\def\tr            {{\mathrm{tr}}}

\def\qed{{\unskip\nobreak\hfil\penalty50
\hskip2em\hbox{}\nobreak\hfil  $\Box$
\parfillskip=0pt \finalhyphendemerits=0\par}\medskip}
\def\proof{\trivlist \item[\hskip \labelsep{\it Proof.\ }]}
\def\endproof{\null\hfill\qed\endtrivlist}

\newcommand\lableq[1]{\label{#1}\end{equation}}
\newcommand\labl[1]{\label{#1}}
\def\typei         {type \nolinebreak I}
\def\typeii        {type \nolinebreak II}
\def\typeiii       {type \nolinebreak III}

\documentclass[11pt]{article}
\usepackage{amssymb,amsfonts,latexsym,epic,eepic}

\begin{document}

%%%%%%%%%%%%%%%%%%%%%%%%%%%%%%%%%%%%%%%%%%%%%%%%%%%%%%%%%%%%%%%%%%%

\newtheorem{definition}{Definition}[section]
\newtheorem{lemma}[definition]{Lemma}
\newtheorem{corollary}[definition]{Corollary}
\newtheorem{theorem}[definition]{Theorem}
\newtheorem{proposition}[definition]{Proposition}
\newtheorem{conjecture}[definition]{Conjecture}
\newtheorem{assumption}[definition]{Assumption}

%%%%%%%%%%%%%%%%%%%%%%%%%%%%%%%%%%%%%%%%%%%%%%%%%%%%%%%%%%%%%%%%%%%%

\title{Modular Invariants from Subfactors:\\
Type I Coupling Matrices and Intermediate Subfactors}
\author{{\sc Jens B\"ockenhauer} and {\sc David E. Evans}\\ \\
School of Mathematics\\University of Wales Cardiff\\
PO Box 926, Senghennydd Road\\Cardiff CF24 4YH, Wales, U.K.}
\maketitle

\begin{abstract}
A braided subfactor determines a coupling matrix $Z$
which commutes with the S- and T-matrices arising from
the braiding. Such a coupling matrix is not necessarily
of ``\typei'', i.e.\ in general it does not have a
block-diagonal structure which can be reinterpreted as
the diagonal coupling matrix with respect to a suitable
extension. We show that there are always two intermediate
subfactors which correspond to left and right maximal
extensions and which determine ``parent'' coupling matrices
$Z^\pm$ of \typei. Moreover it is shown that if the
intermediate subfactors coincide, so that $Z^+=Z^-$,
then $Z$ is related to $Z^+$ by an automorphism of the
extended fusion rules. The intertwining relations of chiral
branching coefficients between original and extended S- and
T-matrices are also clarified. None of our results depends
on non-degeneracy of the braiding, i.e.\ the S- and
T-matrices need not be modular. Examples from $\SOn$ current
algebra models illustrate that the parents can be different,
$Z^+\neq Z^-$, and that $Z$ need not be related
to a \typei\ invariant by such an automorphism.
\end{abstract}

%\newpage
%\tableofcontents

%%%%%%%%%%%%%%%%%%%%%%%%%%%%%%%%%%%%%%%%%%%%%%%%%%%%%%%%%%%%%%%%%%%%%%%%%%%

%%%%%%%%%%%%%%%%%%%%%%%%%%%%%%%%%%%%%%%%%%%%%%%%%%%%%%%%%%%%%%%%%%%%%%%%%%%

\section{Introduction}

A prominent problem in rational conformal field theory
(RCFT) is the classification of modular invariants.
Though it is usually a difficult task and solved only
for a few special models (e.g.\ \cite{CIZ,G2}),
its mathematical formulation is simple: For a given
unitary, finite-dimensional representation of
the modular group $\SLZ$, let $S=(S_{\la,\mu})$
and $T=(T_{\la,\mu})$ denote the matrices representing
the generators ${0\,-1\choose 1\,\,\,0}$ and
${1\,1\choose 0\,1}$, respectively. In the
representations of interest $T$ is diagonal,
$S$ is symmetric, $S^2$ is a permutation,
and $S_{\la,0}\ge S_{0,0}>0$.
Here ``$0$'' is a distinguished label,
referring to the ``vacuum sector''. A modular invariant
is then a ``coupling matrix'' $Z$ (or ``mass matrix'')
commuting with $S$ and $T$,
\begin{equation}
SZ=ZS \,, \qquad \mbox{and} \qquad TZ=ZT \,,
\lableq{commut}
and subject to the constraints
\begin{equation}
Z_{\la,\mu}=0,1,2,\ldots \,, \qquad \mbox{and} \qquad Z_{0,0}=1 \,.
\lableq{constr}
These constraints reflect the physical background
of the problem: The coupling matrix usually describes
multiplicities in the decomposition of the Hilbert space
$\cH_{\mathrm{phys}}$ of physical states of a
2D conformal field theory under the action of a
``symmetry algebra'' $\cA \otimes \cA$ which is
a tensor product of two copies\footnote{In generic
RCFTs one does not necessarily start with a symmetry
algebra made of two identical chiral algebras.
However, \erf{commut} is designed for such a symmetric
situation, whereas in the ``heterotic'' situation one has
to deal with different S- and T-matrices intertwined by Z.}
of a ``chiral algebra'' $\cA$,
\[ \cH_{\mathrm{phys}} = \bigoplus\nolimits_{\la,\mu}
Z_{\la,\mu} \, H_\la \otimes H_\mu \,, \]
giving rise to a modular invariant partition function
\[ \cZ = \sum\nolimits_{\la,\mu}
Z_{\la,\mu} \, \chi_\la \, \chi_\mu^* \,. \]
Here the $\chi_\la$'s and $\chi_\mu$'s are the conformal
characters of the representations $H_\la$ and $H_\mu$,
and the modular group action of $S$ and $T$ comes from
resubstitution of their arguments, leaving the
sesqui-linear combination $\cZ$ invariant.
The condition $Z_{0,0}=1$ then expresses the uniqueness
of the vacuum state.

The simplest example for a coupling matrix is the
``diagonal case'', $Z_{\la,\mu}=\del\la\mu$, which always
gives a modular invariant partition function.
More interesting non-diagonal modular invariants
arise whenever the chiral algebra can be extended
by some local fields. Of special relevance are
the so-called \typei\ invariants \cite{DZ2} for which
the entries of the coupling matrices can be written as
\begin{equation}
Z_{\la,\mu}= \sum\nolimits_\tau b_{\tau,\la} b_{\tau,\mu}\,,
\lableq{ty1}
and which refer directly to the extension
through their ``block-dia\-go\-nal'' structure:
The label $\tau$ runs over the representations of
the extended chiral algebra $\cA^\ext$,
and the non-negative integers $b_{\tau,\la}$ describe
the branching of a representation $\tau$ into $\la$'s
according to the inclusion $\cA\subset\cA^\ext$.
The branching coefficients fulfill $b_{\tau,0}=\del \tau 0$
(by some abuse of notation we denote the vacuum sector of
$\cA^\ext$ also by ``0''),
thus guaranteeing the normalization condition $Z_{0,0}=1$.
Rewriting the partition function in terms of the
extended characters
$\chi^\ext_\tau=\sum_\la b_{\tau,\la} \chi_\la$,
any \typei\ modular invariant can be considered
as completely diagonal:
$Z_{\tau,\tau'}^\ext=\del {\tau'}\tau$.
It is argued in \cite{DV,MS} that after extending the
chiral algebras maximally, the coupling matrix of a
partition function in RCFT is at most a permutation,
$Z_{\tau,\tau'}^\ext=\del {\tau'}{\om(\tau)}$,
where the permutation $\om$ is an automorphism of the
extended fusion rules, satisfying $\om(0)=0$.
As a consequence, a maximal extension
$\cA\subset\cA^\ext$ in RCFT produces
a coupling matrix of some modular invariant
partition function which can be written as
\begin{equation}
Z_{\la,\mu}= \sum\nolimits_\tau b_{\tau,\la} b_{\om(\tau),\mu}\,.
\lableq{ty2}
Partition functions of the form \erf{ty2} which are not
of \typei, \erf{ty1}, are usually referred to as being
``\typeii'' \cite{DZ2}.

Given matrices $S$ and $T$ arising from the modular
transformations of a collection of characters $\chi_\la$
in a RCFT, the solution of the mathematical problem given in
Eqs.\ (\ref{commut}) and (\ref{constr}) can neverthelss yield
coupling matrices which are neither of \typei\ nor of \typeii.
Note that a coupling matrix of the form in \erf{ty2}
has necessarily ``symmetric vacuum coupling'',
$Z_{\la,0}=Z_{0,\la}$. However, even for rather well-behaved
models like $\SOn$ current algebras there are known matrices
$Z$ satisfying Eqs.\ (\ref{commut}) and (\ref{constr}), but which
do not have this symmetry, cf.\ Section \ref{heterotic} below.
Chiral algebras often admit different extensions,
and only then, but much more rarely, modular invariants
without symmetric vacuum coupling have been found.
Namely, it can happen that two chiral extensions
$\cA\subset\cA^\ext_\pm$ of the original chiral algebra
$\cA$ are compatible such that a given coupling matrix has
to be interpreted as an ``automorphism'' invariant with respect
to the enhanced ``heterotic'' symmetry algebra
$\cA^\ext_+\otimes\cA^\ext_-$.
(It seems that this possibility has sometimes been ignored
in the literature although the heterotic case was taken
into account in \cite{MS}.) Unfortunately the standard
terminology ``permutation'' and ``automorphism'' is a bit
misleading in the heterotic case because the labels of left
and right sectors are generically different. A more precise
notion would be ``bijection'' and
``isomorphism of fusion rules'', and the distinction between
diagonal and permutation invariant does no longer make sense
for a maximally extended heterotic symmetry algebra.
Finally, in case that for a fixed theory there are
several modular invariant partition functions it may
happen that a linear combination of their coupling
matrices yields a solution of \erf{constr},
which may however fail to have a consistent
interpretation as a partition function \cite{SY,V,FSS}.
Such modular invariants without physical interpretation
seem to be extremely rare.

The mathematical classification problem of
Eqs.\ (\ref{commut}) and (\ref{constr}) was considered
in \cite{BEK1,BEK2} by means of subfactor theory,
using the ideas of $\a$-induction \cite{LR,X,BE1,BE2,BE3}
and double triangle algebras \cite{O}.
The analysis in \cite{BEK1,BEK2} addressed in particular
the problem of understanding the relation between modular
invariants, graphs and ``nimreps'' (non-negative integer
valued matrix representations of the Verlinde fusion algebra)
--- a puzzling connection going back to the celebrated
A-D-E classification of \cite{CIZ,Kt},
its general nature noticed in \cite{DZ1,DZ2}
and further studied in \cite{DiF,PZ,BPPZ}.
It follows from \cite{R1,FG1,FRS2} that a
(\typeiii) von Neumann factor $N$ with a system $\NXN$
of braided endomorphisms give rise to certain
``statistics'' matrices $S$ and $T$, which are modular
whenever the braiding is non-degenerate \cite{R1}.
It was shown in \cite{BEK1} that then an inclusion
$N\subset M$ of von Neumann factors which is
compatible with the system $\NXN$ determines
a coupling matrix $Z$ by $\a$-induction,
\begin{equation}
Z_{\la,\mu} = \lan \a^+_\la,\a^-_\mu \ran \,, \qquad
\la,\mu\in \NXN \,,
\lableq{Za+a-}
solving Eqs.\ (\ref{commut}) and (\ref{constr})
even if the braiding is degenerate.
Here $\a^+_\la$ and $\a^-_\mu$ are the two inductions
of $\la$ and $\mu$, coming from braiding and opposite
braiding, and the bracket
$\lan \a^+_\la,\a^-_\mu \ran$ denotes the dimension
of their relative intertwiner space. From
current algebra models (``WZW'') in RCFT one can
construct braided subfactors such that the statistics
matrices $S$ and $T$ and the Kac-Peterson matrices
performing the $\SLZ$ transformations of the affine
characters (cf.\ \cite{K}) coincide.
This connection between statistics and conformal character
transformations is expected to hold quite generally in RCFT
(e.g.\ it was conjectured in \cite{FG2}), and the conformal
spin-statistics theorem \cite{FG1,FRS2,GL} establishes
this for the T-matrices. To prove this for the S-matrices
requires one to show that the composition of superselection
sectors indeed recovers the Verlinde fusion rules, and this
has been done for several models, most significantly for
$\SUn$ at all levels in \cite{W}.
For local extensions (cf.\ \cite{RST}) the subfactor
$N\subset M$ can be thought of as a version of the
inclusion $\cA\subset\cA^\ext$.
In terms of (a variation of) the $\alpha$-induction
formula of \cite{LR}, such subfactors were first
investigated for certain conformal inclusions in \cite{X}.
This was further analyzed and extended to simple current
extensions in \cite{BE1,BE2}, and that $\a$-induction
indeed recovers the corresponding modular invariants was
found (for $\SUz$ and $\SUd$ current algebras) in \cite{BE3}.

The two inductions, $\a^+$ and $\a^-$, produce chiral
systems $\MXMp$ and $\MXMm$, intersecting on the
``ambichiral'' system $\MXMo$. Then \erf{Za+a-} can
be written as
\begin{equation}
Z_{\la,\mu}= \sum\nolimits_{\tau\in\MXMo}
b^+_{\tau,\la} b^-_{\tau,\mu}\,,
\lableq{Zb+b-}
with {\em chiral} branching coefficients
$b^\pm_{\tau,\la}=\lan\tau,\a^\pm_\la\ran$. Now the
question arises whether the general subfactor setting
of \cite{BEK1,BEK2}, which is also able to produce
\typeii\ invariants so that in particular
$b^+_{\tau,\la}\neq b^-_{\tau,\la}$ is possible,
will be confined to coupling matrices of the form
of \erf{ty2} or whether it can even produce other
solutions of Eqs.\ (\ref{commut}) and (\ref{constr}),
e.g.\ with heterotic vacuum coupling.
This is the issue of the present paper.
We will indeed demonstrate that our framework incorporates
the general situation, including modular invariants
corresponding to heterotic extensions of the symmetry algebra.

In fact, we study subfactors $N\subset M$, producing
coupling matrices $Z$, through intermediate subfactors,
making essential use of the Galois correspondence elaborated
in \cite{ILP}. We derive in Section \ref{interm} that
there are intermediate subfactors $M_+$ and $M_-$,
$N\subset M_\pm\subset M$, naturally associated to the vacuum
column $(Z_{\la,0})$ respectively the vacuum row $(Z_{0,\la})$
of the coupling matrix determined by $N\subset M$.
The subfactors $N\subset M_\pm$ in turn determine
coupling matrices $Z^\pm$ which are of the form of
\erf{ty1} and can be interpreted as the ``\typei\ parents''
of the original coupling matrix $Z$.
In Section \ref{cfra} we show that in the case $M_+=M_-$,
so that in particular there is a unique parent $Z^+=Z^-$,
the coupling matrix is indeed of the form of \erf{ty2},
recovering a fusion rule automorphism of the ambichiral system.
For the general situation we prove a proposition which
shows that $M_+$ and $M_-$ should be regarded as the operator
algebraic version of maximally extended left and right
chiral algebras, using a recent result of Rehren \cite{R7}.
In Section \ref{SText} we establish the intertwining relations
of the chiral branching coefficients between the original
S- and T-matrices and the ``extended'' ones arising from
the ambichiral braiding.
It is remarkable that the entire analysis does not
need to assume that the braiding is non-degenrate, i.e.\
our results remain valid even if the matrices $S$ and $T$
are not modular.
In Section \ref{heterotic} we finally show by examples from
$\SOn$ current algebras that indeed $M_+\neq M_-$ is possible,
that the parents can be different, $Z^+\neq Z^-$,
and that subfactors can produce coupling matrices $Z$
which have heterotic vacuum coupling,
so that \erf{ty2} can not be adopted in general.

\section{Preliminaries}
\labl{prelim}

Let $A$ and $B$ be \typeiii\ von Neumann factors.
A unital $\ast$-homomorphism $\rho:A\rightarrow B$
is called a $B$-$A$ morphism. The positive number
$d_\rho=[B:\rho(A)]^{1/2}$ is called the statistical
dimension of $\rho$; here $[B:\rho(A)]$ is the
Jones index \cite{J1} of the subfactor
$\rho(A)\subset B$. If $\rho$ and $\sig$ are $B$-$A$
morphisms with finite statistical dimensions, then
the vector space of intertwiners
\[ \Hom(\rho,\sig)=\{ t\in B: t\rho(a)=\sig(a)t \,,
\,\, a\in A \}  \]
is finite-dimensional, and we denote its dimension by
$\lan\rho,\sig\ran$.
An $A$-$B$ morphism $\co\rho$ is a
conjugate morphism if there are isometries
$r_\rho\in\Hom(\id_A,\co\rho\rho)$ and
${\co r}_\rho\in\Hom(\id_B,\rho\co\rho)$ such that
$\rho(r_\rho)^* {\co r}_\rho=d_\rho^{-1}\bfe_B$ and
$\co\rho({\co r}_\rho)^* r_\rho=d_\rho^{-1}\bfe_A$.
The map $\phi_\rho:B\rightarrow A$,
$b\mapsto r_\rho^* \co\rho(b)r_\rho$, is called the
(unique) standard left inverse and satisfies
\be
\phi_\rho(\rho(a)b\rho(a'))=a\phi_\rho(b)a' \,,
\quad a,a'\in A\,,\quad b\in B \,.
\lableq{phiaba}
If $t\in\Hom(\rho,\sig)$ then we have
\be
d_\rho \phi_\rho (bt) = d_\sig \phi_\sig (tb) \,,
\qquad b\in B \,.
\lableq{phibttb}

We work with the setting of \cite{BEK1}, i.e.\ we are
working with a \typeiii\ subfactor and
finite system $\NXN\subset\End(N)$ of (possibly degenerately)
braided morphisms which is compatible with the inclusion
$N\subset M$. Then the inclusion is in particular forced to have
finite Jones index and also finite depth (see e.g.\ \cite{EK}).
More precisely, we make the following

\begin{assumption}
{\rm We assume that we have a
\typeiii\ subfactor $N\subset M$
together with a finite system of endomorphisms
$\NXN\subset\End(N)$ in the sense of \cite[Def.\ 2.1]{BEK1}
which is braided in the sense of \cite[Def.\ 2.2]{BEK1}
and such that $\canr=\co\iota\iota\in\Sigma(\NXN)$ for the
injection $M$-$N$ morphism $\iota:N\hookrightarrow M$ and a
conjugate $N$-$M$ morphism $\co\iota$.}
\labl{assbasic}
\end{assumption}

With the braiding $\eps$ on $\NXN$ and its
extension to $\Sigma(\NXN)$ (the set of finite sums of
morphisms in $\NXN$) as in \cite{BEK1}, one can
define the $\a$-induced morphisms $\a^\pm_\la\in\End(M)$
for $\la\in\Sigma(\NXN)$ by the Longo-Rehren formula \cite{LR},
namely by putting
\[ \a_\la^\pm = \co\iota^{\,-1} \circ \Ad
(\eps^\pm(\lambda,\canr)) \circ \lambda \circ \co\iota \,, \]
where $\co\iota$ denotes a conjugate morphism of the
injection map $\iota:N\hookrightarrow M$.
Then $\a^+_\la$ and $\a^-_\la$ extend $\la$, i.e.\
$\a^\pm_\la\circ\iota=\iota\circ\la$, which in turn implies
$d_{\a_\la^\pm}=d_\la$ by the multiplicativity of
the minimal index \cite{L3}. 
Let $\can=\iota\co\iota$ denote Longo's canonical
endomorphism from $M$ into $N$. Then there is an isometry
$v\in\Hom(\id,\can)$ such that any $m\in M$ is uniquely
decomposed as $m=nv$ with $n\in N$.
Thus the action of the extensions
$\a^\pm_\la$ are uniquely characterized by the relation
$\a^\pm_\la(v)=\eps^\pm(\la,\canr)^* v$ which can be
derived from the braiding fusion equations
(BFE's, see e.g.\ \cite[Eq.\ (5)]{BEK1}).
Moreover, we have $\a_{\la\mu}^\pm=\a_\la^\pm \a_\mu^\pm$
if also $\mu\in\Sigma(\NXN)$, and clearly
$\a_{{\rm{id}}_N}^\pm={{\rm{id}}}_M$.
In general one has
\[ \Hom(\la,\mu) \subset \Hom(\a^\pm_\la,\a^\pm_\mu)
\subset \Hom(\iota\la,\iota\mu) \,,\qquad\la,\mu\in\Sigma(\NXN)\,. \]
The first inclusion is a consequence of the BFE's.
Namely, $t\in\Hom(\la,\mu)$ obeys
$t\eps^\pm(\canr,\la)=\eps^\pm(\canr,\mu)\canr(t)$,
and thus
\[ t \a_\la^\pm(v) = t \eps^\pm(\la,\canr)^* v =
\eps^\pm(\mu,\canr)^* \canr(t) v
=\eps^\pm(\mu,\canr)^* vt = \a_\mu^\pm(v) t \,. \]
The second follows from the extension property
of $\a$-induction. Hence $\a_{\co\la}^\pm$
is a conjugate for $\a_\la^\pm$ as there are
$r_\la\in\Hom(\id,\co\la\la)\subset
\Hom(\id,\a_{\co\la}^\pm\a_\la^\pm)$
and ${\co r}_\la\in\Hom(\id,\la\co\la)\subset
\Hom(\id,\a_\la^\pm\a_{\co\la}^\pm)$ such
that
$\la(r_\la)^*{\co r}_\la=\co\la({\co r}_\la)^*r_\la
= d_\la^{-1}\bfe$. We also have some kind of
naturality equations for $\a$-induced morphisms,
\begin{equation}
x \eps^\pm(\rho,\la)=\eps^\pm(\rho,\mu)\a^\pm_\rho(x)
\lableq{anat}
whenever $x\in\Hom(\iota\la,\iota\mu)$
and $\rho\in\Sigma(\NXN)$.

Recall that the statistics phase of $\om_\la$ for
$\la\in\NXN$ is given as
\[ d_\la \phi_\la(\eps^+(\la,\la))=\om_\la \bfe \,. \]
The monodromy matrix $Y$ is defined by
\[ Y_{\la,\mu} = \sum_{\rho\in\NXN}
\frac{\om_\la \om_\mu}{\om_\rho} N_{\la,\mu}^\rho d_\rho \,,
\qquad \la,\mu\in\NXN \,,\]
with $N_{\la,\mu}^\rho=\lan\rho,\la\mu\ran$ denoting the
fusion coefficients. Then one checks that $Y$ is symmetric,
that $Y_{\co\la,\mu}=Y_{\la,\mu}^*$ as well as
$Y_{\la,0}=d_\la$ \cite{R1,FG1,FRS2}. (As usual, the
label ``$0$'' refers to the identity morphism $\id\in\NXN$.)
Now let $\Omega$ be the diagonal matrix with entries
$\Omega_{\la,\mu}=\om_\la \del\la\mu$.
Putting $Z_{\la,\mu}=\lan\a^+_\la,\a^-_\mu\ran$
defines a matrix subject to the constraints \erf{constr}
and commuting with $Y$ and $\Omega$ \cite{BEK1}.
The Y- and $\Omega$-matrices obey
$\Omega Y \Omega Y \Omega = z Y$
where $z=\sum_\la d_\la^2 \om_\la$ \cite{R1,FG1,FRS2},
and this actually holds even if the braiding is
degenerate (see \cite[Sect.\ 2]{BEK1}).
If $z\neq 0$ we put $c=4\arg(z)/\pi$, which is defined
modulo 8, and call it the ``central charge''.
Moreover, S- and T-matrices are then defined by
\[ S = |z|^{-1} Y \,, \qquad T = \E^{-\I\pi c/12} \Omega \]
and hence fulfill $TSTST=S$. 
One has $|z|^2=w$ with the global index
$w=\sum_\la d_\la^2$ and $S$ is unitary, so that
$S$ and $T$ are indeed the standard generators in a
unitary representation of the modular group $\SLZ$,
if and only if the braiding is non-degenerate \cite{R1}.
Consequently, $Z$ gives a modular invariant in this case.

Let $\MXM\subset\End(M)$ denote a system of endomorphisms
consisting of a choice of representative endomorphisms of
each irreducible subsector of sectors of the form
$[\iota\la\co\iota]$, $\la\in\NXN$.
We choose $\id\in\End(M)$
representing the trivial sector in $\MXM$.
Then we define similarly the chiral systems
$\MXMpm$ and the $\a$-system $\MXMa$ to be the subsystems
of endomorphisms $\beta\in\MXM$
such that $[\beta]$ is a subsector of $[\a^\pm_\la]$ and of
of $[\a_\la^+\a_\mu^-]$, respectively,
for some $\la,\mu\in\NXN$.
(Note that any subsector of $[\a_\la^+\a_\mu^-]$ is
automatically a subsector of  $[\iota\nu\co\iota]$
for some $\nu\in\NXN$.)
The ambichiral system is defined
as the intersection $\MXMo=\MXMp\cap\MXMm$, so that
$\MXMo \subset \MXMpm \subset \MXMa \subset \MXM$.
Their ``global indices'', i.e.\ their sums over the
squares of the statistical dimensions are denoted
by $w_0$, $w_\pm$, $w_\a$ and $w$, and thus fulfill
$1 \le w_0 \le w_\pm \le w_\a \le w$.

\section{More on global indices and chiral locality}
\labl{globinloc}

Recall from \cite[Prop.\ 3.1]{BEK2} that
$w_+=w_-$ and that
$w/w_+=\sum_{\la\in\NXN} d_\la Z_{\la,0}$.
We will now derive a general formula for the $\a$-global index
$w_\a=\sum_{\beta\in\MXMa}d_\beta^2$ and also for $w_0$.
We denote by $\NXNd\subset\NXN$ the subsystem of
degenerate morphisms.

\begin{proposition}
The $\a$-global index is given by
\be
w_\a = \frac w{\sum_{\la\in\NXNd} Z_{0,\la} d_\la} .
\lableq{wa}
Moreover, the ambichiral global index is given by
$w_0=w_+^2/w_\a$.
\labl{aglob}
\end{proposition}

\proof
Let $R_{\la,\mu}$, $\la,\mu\in\NXN$, denote matrices
with entries
$R_{\la,\mu;\beta}^{\beta'}=
\lan \beta\a_\la^+\a_\mu^-,\beta'\ran$,
$\beta,\beta'\in\MXMa$. Further let $\vec{d}$ denote
the column vector with entries $d_\beta$, $\beta\in\MXMa$.
Then $\vec{d}$ is a simultaneous eigenvector of the
matrices $R_{\la,\mu}$ with respective eigenvalues
$d_\la d_\mu$. We define another vector $\vec{v}$
by putting
\[ v_\beta = \sum_{\la,\mu\in\NXN} d_\la d_\mu
\lan \beta, \a^+_\la \a^-_\mu \ran \,,
\qquad \beta\in\MXMa \,. \]
Then we have $R_{\la,\mu}\vec{v}=d_\la d_\mu\vec{v}$,
as we can compute
\[ \begin{array}{ll}
(R_{\la,\mu} \vec{v})_\beta &= \sum_{\beta'\in\MXMa}
\sum_{\nu,\rho\in\NXN} \lan \beta \a_\la^+\a_\mu^-,
\beta' \ran d_\nu d_\rho \lan \beta',
\a_\nu^+ \a_\rho^- \ran  \\[.4em]
&= \sum_{\nu,\rho\in\NXN} d_\nu d_\rho
\lan \beta \a_\la^+\a_\mu^-,\a_\nu^+ \a_\rho^- \ran \\[.4em]
&=  \sum_{\nu,\rho,\xi,\eta\in\NXN} d_\nu d_\rho
N_{\nu,\co\la}^\xi N_{\rho,\co\mu}^\eta
\lan \beta , \a_\xi^+ \a_\eta^- \ran  \\[.4em]
&=  \sum_{\xi,\eta\in\NXN} d_\la d_\mu d_\xi d_\eta
\lan \beta , \a_\xi^+ \a_\eta^- \ran
= d_\la d_\mu v_\beta \,.
\end{array} \]
Because the sum matrix $\sum_{\la,\mu} R_{\la,\mu}$ is
irreducible it follows $\vec{v}=\zeta \vec{d}$,
$\zeta\in\bbR$, by the uniqueness of the Perron-Frobenius
eigenvector. Note that
$d_\la d_\mu=\sum_\beta\lan\beta,\a_\la^+\a_\mu^-\ran d_\beta$,
and hence
$w^2=\sum_\beta v_\beta d_\beta = \zeta w_\a$.
We next notice that $\zeta=v_0$ as $d_0=1$. But $v_0$
can be computed as
\[ v_0 = \sum_{\la,\mu\in\NXN} d_\la d_\mu
\lan \a_{\co\la}^+, \a_\mu^-\ran =
\sum_{\la,\mu\in\NXN} Y_{0,\la} Z_{\la,\mu} Y_{\mu,0}
= \sum_{\la,\mu\in\NXN} Z_{0,\la} Y_{\la,\mu} Y_{\mu,0} \,,\]
where we used commutativity of the monodromy matrix $Y$
with the coupling matrix $Z$ \cite[Thm.\ 5.7]{BEK1}.
By Rehren's argument \cite{R1} we have
\[ \sum_{\mu\in\NXN} Y_{\la,\mu} Y_{\mu,0} = \left\{
\begin{array}{cl} w d_\la & \qquad \la\in\NXNd \\
0 & \qquad \la\notin\NXNd \end{array} \right. . \]
Hence $\zeta=w\sum_{\la\in\NXNd} Z_{0,\la}d_\la$,
establishing \erf{wa}.

Next we define two vectors $\vec{v}^\pm$ with entries
$v^\pm_\beta = \sum_\la d_\la \lan\beta,\a^\pm_\la\ran$,
$\beta\in\MXMa$. From
(the proof of) \cite[Prop.\ 3.1]{BEK2} we learn that
\[ v^\pm_\beta =  \left\{ \begin{array}{cl}  d_\beta w / w_+
& \qquad \beta\in\MXMpm \\ 0 & \qquad \beta\notin\MXMpm
\end{array} \right. . \]
Consequently
$\lan \vec{v}^+,\vec{v}^- \ran = w_0 w^2/w_+^2$.
But we can also compute directly
\[ \lan \vec{v}^+,\vec{v}^- \ran 
=\sum_{\la,\mu} \sum_{\beta\in\MXMo} d_\la
\lan \a^+_\la,\beta \ran \lan \beta, \a^-_\mu \ran d_\mu
= \sum_{\la,\mu}
d_\la Z_{\la,\mu} d_\mu = \zeta = \frac{w^2}{w_\a} \,.\]
completing the proof.
\endproof

Note that Proposition \ref{aglob} in particular provides
a new proof of the ``generating property of $\a$-induction'',
i.e.\ $\MXMa=\MXM$ if the braiding is non-degenerate,
which was established in \cite[Thm.\ 5.10]{BEK1}.

Now recall that the chiral locality condition
$\eps^+(\canr,\canr)v^2=v^2$ expresses local commutativity
(``locality'') of the extended net, if $N\subset M$
arises from a net of subfactors \cite{LR}.

\begin{proposition}
The following conditions are equivalent:
\begin{enumerate}
\item We have $Z_{\la,0}=\lan\canr,\la\ran$ for all $\la\in\NXN$.
\item We have $Z_{0,\la}=\lan\canr,\la\ran$ for all $\la\in\NXN$.
\item Chiral locality holds: $\eps^+(\canr,\canr)v^2=v^2$.
\end{enumerate}
\labl{condloc}
\end{proposition}

\proof
The implications {\it 3} $\Rightarrow$ {\it 1,2}
follow from \cite[Thm.\ 3.9]{BE1}. We need to show
{\it 1,2} $\Rightarrow$ {\it 3}. Recall
$\lan\canr,\la\ran=\lan\iota,\iota\la\ran$. Moreover,
by the extension property of $\a$-induction we have
\[ \Hom(\id,\a^\pm_\la) \subset \Hom(\iota,\iota\la) \,,
\qquad \la\in\NXN \,. \]
Hence, if $Z_{\la,0}=\lan\canr,\la\ran$
(respectively $Z_{0,\la}=\lan\canr,\la\ran$) then
$\Hom(\id,\a^\pm_\la) = \Hom(\iota,\iota\la)$ for all
$\la\in\NXN$. Then consequently
$\Hom(\id,\a^\pm_\canr) = \Hom(\iota,\iota\canr)$.
We clearly have $v\in\Hom(\iota,\iota\canr)$ and hence
$v^2 = \a^\pm_\canr(v)v = \eps^\pm(\canr,\canr)^* v^2$,
i.e.\ chiral locality holds.
\endproof

Recall from \cite[Prop.\ 3.4]{BEK2} that the coupling matrix
arising from a braided subfactor with satisfied chiral
locality condition is automatically of \typei. Hence
Proposition \ref{condloc} states that chiral locality
is equivalent to the canonical endomorphism being
``fully visible'' in the vacuum row (or column) of the
coupling matrix.

\section{Intermediate subfactors}
\labl{interm}

In this section we are searching for certain intermediate
subfactors
\[ N\subset \tilde M \subset M \]
of our subfactor $N\subset M$. It follows from
\cite[Sect.\ 3]{ILP} that the set of such intermediate
subfactors $\tilde M$ is in a bijective correspondence
with systems of subspaces $\cK_\rho\subset\cH_\rho$,
where $\cH_\rho=\Hom(\iota,\iota\rho)$, $\rho\in\NXN$,
and subject to conditions
\begin{enumerate}
\item[(i.)] $\cK_\rho^* \subset N \cK_{\co\rho}$,
\item[(ii.)] $\cK_\rho\cK_\sig \subset
 \sum_{\xi\prec\rho\sig} N \cK_\xi$,
\end{enumerate}
where the sum in (ii.) runs over all $\xi\in\NXN$
such that $N_{\rho,\sig}^\xi >0$. The factor
$\tilde M$ is then generated by $N$ and the $\cK_\rho$'s
and is uniquely decomposed as
\[ \tilde M = \sum_\rho N \cK_\rho \,. \]
The dual canonical endomorphism $\tilde \canr$ of
$N\subset \tilde M$ decomposes as a sector as
$[\tilde \canr] = \bigoplus_\rho n_\rho [\rho]$,
where $n_\rho=\dim\cK_\rho$.

We now define the spaces
\[ \cK_\rho^\pm = \Hom(\id,\a_\rho^\pm) \,, \qquad
\rho\in\NXN \,. \]
Note that
$\cK_\rho^\pm\subset\cH_\rho=\Hom(\iota,\iota\rho)$,
that $\dim\cK^+_\rho=Z_{\rho,0}$ and
$\dim\cK^-_\rho=Z_{0,\rho}$.

\begin{lemma}
We have
\begin{equation}
\cK_\rho^\pm = \{ zv : z\in\Hom(\canr,\rho) \,,\,\,
z \can(v) = z \eps^\mp(\canr,\canr) \can(v) \}
\end{equation}
\labl{cK=zv}
\end{lemma}

\proof
Let $x\in\cK_\rho^\pm$. Now $x$ is uniquely decomposed
as $x=zv$ with $z\in N$. Clearly $z\in\Hom(\canr,\rho)$.
Then $x\in\Hom(\id,\a^\pm_\rho)$ reads,
using naturality (see e.g.\ \cite[Eq.\ (8)]{BEK1}),
\[ z\can(v)v = zv^2 = \a^\pm_\rho(v)zv
= \eps^\mp(\canr,\rho) vzv
= \eps^\mp(\canr,\rho)\canr(z)\can(v)v
= z \eps^\mp(\canr,\canr)\can(v)v \,, \]
hence is equivalent to
$z\can(v)=z \eps^\mp(\canr,\canr)\can(v)$.
\endproof

We choose orthonormal bases of isometries
$t(_{\rho,\sig}^\xi)_i \in\Hom(\xi,\rho\sig)$,
$i=1,2,...,N_{\rho,\sig}^\xi$ so that
$\sum_\xi \sum_{i=1}^{N_{\rho,\sig}^\xi}
t(_{\rho,\sig}^\xi)_i t(_{\rho,\sig}^\xi)_i^*=\bfe$.

\begin{lemma}
We have
\begin{enumerate}
\item[{\rm (i.)}] $(\cK_\rho^\pm)^* = r_\rho^* \cK_{\co\rho}^\pm$,
\item[{\rm (ii.)}] $\cK_\rho^\pm \cK_\sig^\pm \subset
 \sum_\xi \sum_{i=1}^{N_{\rho,\sig}^\xi}
 t(_{\rho,\sig}^\xi)_i \cK_\xi^\pm$.
\end{enumerate}
\labl{KKK}
\end{lemma}

\proof
Right Frobenius reciprocity \cite{I2} gives us an isomorphism
$\Hom(\id,\a^\pm_{\co\rho})\rightarrow\Hom(\id,\a^\pm_\rho)$,
$x\mapsto x^*r_\rho$. Thus any element
$y\in\cK_\rho^\pm$ can be written as $x^*r_\rho$
with some $x\in\cK_{\co\rho}^\pm$,
proving (i.). To prove (ii.), let
$x_\rho=z_\rho v\in\cK_\rho^\pm$ and
$x_\sig = z_\sig v\in\cK_\sig^\pm$ be the decompositions
according to Lemma \ref{cK=zv}. Then
\[ x_\rho x_\sig = z_\rho v z_\sig v =
\sum_\xi \sum_{i=1}^{N_{\rho,\sig}^\xi}
t(_{\rho,\sig}^\xi)_i  t(_{\rho,\sig}^\xi)_i^*
z_\rho \canr(z_\sig) \can(v)v \,. \]
We first notice that
$t(_{\rho,\sig}^\xi)_i^* z_\rho \canr(z_\sig)
\can(v)\in\Hom(\canr,\xi)$.
Next we check by use of the BFE and by $v^2=\can(v)v$ that
\[ \begin{array}{ll}
z_\rho \canr(z_\sig) \can(v) \cdot
\eps^\mp(\canr,\canr) \can(v)
&= z_\rho \canr(z_\sig) \eps^\mp(\canr,\canr^2)
\canr(\can(v))\can(v) \\[.4em]
&= z_\rho \canr(z_\sig) \canr(\eps^\mp(\canr,\canr))
\eps^\mp(\canr,\canr)\can(v)^2 \\[.4em]
&=\rho(z_\sig \eps^\mp(\canr,\canr))
 z_\rho \eps^\mp(\canr,\canr)\can(v)^2 \\[.4em]
&= z_\rho \canr(z_\sig) \canr(\eps^\mp(\canr,\canr))
\can(v)^2\\[.4em]
&= z_\rho \canr(z_\sig\eps^\mp(\canr,\canr)
\can(v))\can(v)
=  z_\rho \canr(z_\sig) \can(v) \cdot \can(v) \,.
\end{array} \]
We conclude
$x_\rho x_\sig= \sum_\xi \sum_{i=1}^{N_{\rho,\sig}^\xi}
t(_{\rho,\sig}^\xi)_i x_\xi^i$ with
$x_\xi^i = t(_{\rho,\sig}^\xi)_i^*
z_\rho \canr(z_\sig) \can(v)v\in\cK_\xi^\pm$.
\endproof

\begin{corollary}
There are two (possibly identical) intermediate
subfactors $N\subset M_\pm\subset M$ with
$M_\pm=\sum_\rho N \cK_\rho^\pm$.
\labl{Mpm}
\end{corollary}

Our next aim is to show that the subfactors
$N\subset M_\pm$ obey the chiral locality condition.
Let $\iota_\pm:N\hookrightarrow M_\pm$ denote the
injection maps, so that the (dual) canonical endomorphisms
are given by $\can_\pm=\iota_\pm{\co\iota}_\pm$ and
$\canr_\pm={\co\iota}_\pm\iota_\pm$.
We now know that
$[\canr_+]=\bigoplus_\rho Z_{\rho,0}[\rho]$ and
$[\canr_-]=\bigoplus_\rho Z_{0,\rho}[\rho]$.
Due to commutativity of $Y$ and $Z$, which yields in
particular
$\sum_\rho d_\rho Z_{\rho,0}=\sum_\rho Z_{0,\rho} d_\rho$,
we find $d_{\canr_+}=d_{\canr_-}$, i.e.\ the subfactors
$N\subset M_+$ and $N\subset M_-$ have the same Jones index.
Moreover, we can apply $\a$-induction, i.e.\ we
define morphisms in $\End(M_\delta)$ by
\[ \alt \delta \pm \la  = \co\iota_\delta^{\,-1} \circ \Ad
(\eps^\pm(\lambda,\canr_\delta))
\circ \lambda \circ \co\iota_\delta \,, \]
where the index $\delta$ is either $\delta=+$ or $\delta=-$.
This will give rise to ``parent'' coupling matrices $Z^\delta$.
Thanks to Prop.\ \ref{condloc}, it suffices to show
$Z^+_{\la,0}= Z_{\la,0}$ and $Z^-_{0,\la}=Z_{0,\la}$,
$\la\in\NXN$, to prove that the chiral locality condition
holds for $N\subset M_+$ and $N\subset M_-$, respectively.

\begin{lemma}
We have
$\alt \delta \pm \la (x_\rho) = \eps^\pm (\la,\rho)^* x_\rho$
for any $x_\rho\in\cK_\rho^\delta$ and $\la,\rho\in\NXN$.
\labl{a(K)}
\end{lemma}

\proof
Let $x_\rho=z_\rho v$ be the decomposition according to
Lemma \ref{cK=zv}.  We first notice that
\[ \can_\delta (z_\rho v) \canr_\delta (n) =
\can_\delta (z_\rho \canr(n) v) =
\can_\delta (\rho(n) z_\rho v) =
\canr_\delta \rho(n) \can_\delta (z_\rho v) \]
for any $n\in N$, i.e.\
$\can_\delta (x_\rho)\in\Hom(\canr_\delta,\canr_\delta \rho)$.
Therefore we can compute
\[ \begin{array}{ll}
\can_\delta (\alt \delta \pm \la (x_\rho))
&= \eps^\pm (\la,\canr_\delta) \la\can_\delta(x_\rho)
\eps^\mp (\canr_\delta,\la)
= \eps^\mp (\canr_\delta,\la)^* \eps^\mp
(\canr_\delta\rho,\la) \can_\delta (x_\rho) \\[.4em]
&=  \eps^\mp (\canr_\delta,\la)^*  \eps^\mp (\canr_\delta,\la)
\canr_\delta (\eps^\mp (\rho,\la)) \can_\delta (x_\rho)
= \can_\delta(\eps^\pm (\la,\rho)^* x_\rho) \,,
\end{array} \]
and application of $\can_\delta^{-1}$ yields the statement.
\endproof

In the same manner we obtain of course also
$\a^\pm_\la(x_\rho)=\eps^\pm (\la,\rho)^* x_\rho$
for any $x_\rho\in\Hom(\iota,\iota\rho)$.
Therefore we obtain immediately

\begin{corollary}
We have $\a^\pm_\la |_{M_\delta} = \alt \delta\pm\la$.
\labl{res}
\end{corollary}

Hence
$x_\la x_\rho=\a^\pm_\la(x_\rho)x_\la=\alt\delta\pm\la(x_\rho)x_\la$
whenever $x_\la\in\cK_\la^\pm$ and $x_\rho\in\cK_\rho^\delta$,
and thus in particular
$\cK_\la^\pm\subset \Hom(\id_{M_\pm},\alt\pm\pm\la)$.
(Warning: Note that the super- and subscripts, referring to the
$\pm$-induction respectively to the choice of the algebra $M_\pm$,
now have to be {\em the same} because we have
$\cK_\la^+\subset M_+$ and  $\cK_\la^-\subset M_-$,
but in general not the other way round. Here and in the
following, any formula containing such combined
$\pm$-indices has to be read in such a way that we either
take all the upper or all the lower signs.)
On the other hand we have
\[ \lan \id_{M_\pm},\alt\pm\pm\la \ran \le
\lan \canr_\pm,\la\ran = \dim \cK_\la^\pm \,,\]
Therefore we arrive in particular at

\begin{corollary}
We have $\Hom(\id_{M_\pm},\alt\pm\pm\la)=\cK_\la^\pm$ for
any $\la\in\NXN$.
\labl{Zpml}
\end{corollary}

Corollary \ref{Zpml} tells us in particular that
$Z^+_{\la,0}=Z_{\la,0}$ and $Z^-_{0,\la}=Z_{0,\la}$,
so that both subfactors $N\subset M_\pm$ must satisfy
the chiral locality condition. In turn they must be
irreducible by the argument of \cite[Cor.\ 3.6]{BE1}.
We summarize the discussion in the following

\begin{theorem}
There are two (possibly identical) intermediate
subfactors $N\subset M_\pm\subset M$.
The irreducible subfactors $N\subset M_+$ and
$N\subset M_-$ have the same Jones index, they both
satisfy the chiral locality condition and consequently
give rise to \typei\ coupling matrices $Z^\pm$.
The latter are related to the coupling matrix of
the full subfactor $N\subset M$ through
$Z^+_{\la,0}=Z_{\la,0}$ and $Z^-_{0,\la}=Z_{0,\la}$.
\labl{type1Zs}
\end{theorem}

\section{Chiral fusion rule automorphisms}
\labl{cfra}

We now investigate the relations between the chiral
systems $\MXMp$, $\MXMm$ and $\MXMpp$, $\MXMmm$,
respectively. Note that by \cite[Prop.\ 3.1]{BEK1}
and Theorem \ref{type1Zs}, they all must have the
same chiral global index $w_+$.

\begin{lemma}
We have
\begin{equation}
\Hom(\a^\pm_\la,\a^\pm_\mu)=\Hom(\alt\pm\pm\la,\alt\pm\pm\mu)\,,
\end{equation}
in particular $\Hom(\a^\pm_\la,\a^\pm_\mu)\subset M_\pm$,
for any $\la,\mu\in\Sigma(\NXN)$.
\labl{Hom=Homt}
\end{lemma}

\proof
First let $\xi\in\Sigma(\NXN)$. Then there are orthonormal
bases of isometries $s_{\nu,i}\in\Hom(\nu,\xi)$, $\nu\in\NXN$,
$i=1,2,...,\lan\nu,\xi\ran$, such that
$\sum_{\nu,i} s_{\nu,i}s_{\nu,i}^*=\bfe$.
We may write $x\in\Hom(\id,\a^\pm_\xi)$ as
$\sum_{\nu,i} s_{\nu,i}s_{\nu,i}^* x$, and we notice
$s_{\nu,i}^* x\in\Hom(\id,\a^\pm_\nu)=\Hom(\id,\alt\pm\pm\nu)$,
thanks to Corollary \ref{Zpml}, so that
$x\in\Hom(\id,\alt\pm\pm\xi)$. The same argument
works vice versa. Now let $\la,\mu\in\Sigma(\NXN)$.
Then we have Frobenius isomorphisms
\[ \Hom(\a^\pm_\la,\a^\pm_\mu)  \longrightarrow 
\Hom(\id,\a^\pm_{\co\la\mu})=\Hom(\id,\alt\pm\pm{\co\la\mu})
 \longrightarrow 
\Hom(\alt\pm\pm\la,\alt\pm\pm\mu) \]
which map
\[ t  \mapsto s=\sqrt{d_\la/d_\mu}
\a^\pm_{\co\la}(t) r_\la \,,  \qquad s  \mapsto 
t'= \sqrt{d_\la d_\mu} {\co r}_\la^*
\alt\pm\pm\la(s) \,. \]
As we have
\[ t'= d_\la {\co r}_\la^* \alt\pm\pm\la
(\a^\pm_{\co\la}(t)r_\la) =
d_\la {\co r}_\la^* \a^\pm_{\la\co\la}(t)\la(r_\la)=t \]
it follows
$\Hom(\a^\pm_\la,\a^\pm_\mu)=
\Hom(\alt\pm\pm\la,\alt\pm\pm\mu)\subset M_\pm$.
\endproof

\begin{lemma}
Each $\beta_\pm\in\MXMpm$ is equivalent to an extension
of some $\tilde{\beta}_\pm\in\MXMpmpm$. This association
gives rise to bijections
$\vartheta_\pm:\MXMpm\rightarrow\MXMpmpm$.
\labl{bij}
\end{lemma}

\proof
Assume $\beta\equiv\beta_+\in\MXMp$, i.e.\ there is a
$\la\in\NXN$ and an isometry $t\in\Hom(\beta,\a^+_\la)$.
Then $tt^*$ is a minimal projection in
$\Hom(\alt ++\la,\alt ++\la)$ by Lemma \ref{Hom=Homt}.
Hence there is a $\tilde\beta\in\MXMpp$ and an
isometry $\tilde t\in\Hom(\tilde\beta,\alt ++\la)$
such that $\tilde t \tilde{t}^*=tt^*$.
Thus putting $\beta'(m)=\tilde{t}^*\a^+_\la(m)\tilde t$
for $m\in M$ gives an equivalent endomorphism,
as $\beta=\Ad(u)\circ\beta'$ with the unitary
$u=t^*\tilde t\in M$, and we clearly have
$\beta'|_{M_+}=\tilde\beta$, thanks to Corollary \ref{res}.
It remains to show that, if
$\tilde{\beta}_1,\tilde{\beta}_2 \in \MXMpp$ correspond this
way to different $\beta_1,\beta_2\in\MXMp$, then
$\tilde{\beta}_1$ and $\tilde{\beta}_1$ are disjoint.
Let $t_j\in\Hom(\beta_j,\a^+_{\la_j})$ and
$\tilde{t}_j\in\Hom(\tilde{\beta}_j,\alt ++{\la_j})$
be isometries as above, i.e.\
$t_jt_j^*=\tilde{t}_j\tilde{t}_j^*$, $j=1,2$.
Assume for contradiction that there is a unitary
$q\in\Hom(\tilde{\beta}_1,\tilde{\beta}_2)$.
But then
$\tilde{t}_2 q \tilde{t}_1^*\in
\Hom(\alt ++{\la_1},\alt ++{\la_2})
=\Hom(\a^+_{\la_1},\a^+_{\la_2})$,
so that $t_2^*\tilde{t}_2 q \tilde{t}_1^*t_1$
is a unitary in $\Hom(\beta_1,\beta_2)$,
in contradiction to $\beta_1,\beta_2$ being
different elements in $\MXMpm$. Hence the
association $\beta\mapsto\tilde\beta$ defines
a bijection
$\vartheta_+:\MXMp\rightarrow\MXMpp$.
The proof is completed by
exchanging ``$+$'' by ``$-$'' signs.
\endproof

\begin{lemma}
The bijections
$\vartheta_\pm:\MXMpm\rightarrow\MXMpmpm$ preserve
the chiral bran\-ching,
\begin{equation}
\lan\beta,\a^\pm_\la\ran = \lan \vartheta_\pm(\beta),
\alt \pm\pm\la \ran \,, \qquad \la\in\NXN \,, \quad
\beta\in\MXMpm \,,
\lableq{cbran}
and the chiral fusion rules
\begin{equation}
\lan\beta_3,\beta_1\beta_2\ran=\lan \vartheta_\pm(\beta_3),
\vartheta_\pm(\beta_1)\vartheta_\pm(\beta_2) \ran \,, \qquad
\beta_1,\beta_2,\beta_3\in\MXMpm \,,
\lableq{cfus}
and the statistical dimensions.
\labl{furuiso}
\end{lemma}

\proof
We just consider the ``$+$'' case,
the proof for ``$-$'' is analogous.
By Lemma \ref{bij} we may and do assume for simplicity that
now all $\beta\in\MXMp$ are choosen such that
$\beta|_{M_+}=\vartheta_+(\beta)$. This already
forces equality of statistical dimensions
$d_\beta=d_{\vartheta_+(\beta)}$.
Moreover, we just have
$\Hom(\beta,\a^+_\la) = \Hom(\vartheta_+(\beta),
\alt ++\la)$, giving \erf{cbran}.
Given isometries
$t_j\in\Hom(\beta_j,\a^+_{\la_j})
=\Hom(\vartheta_+(\beta_j),\alt ++{\la_j})$,
$j=1,2,3$, and also
$y\in\Hom(\beta_3,\beta_1\beta_2)$, then we find similarly
\[ t_1 \a^+_{\la_1}(t_2)yt_3^* \in
\Hom(\a^+_{\la_3},\a^+_{\la_1}\a^+_{\la_2})=
\Hom(\alt ++{\la_3},\alt ++{\la_1}\alt ++{\la_2}) \]
by Lemma \ref{Hom=Homt}. Therefore, by using
$\a^+_{\la_1}(t_2)=\alt ++{\la_1}(t_2)$
we finally find that $y\in\Hom(\vartheta_+(\beta_3),
\vartheta_+(\beta_1)\vartheta_+(\beta_2))$.
The same argument works vice versa, so that the
intertwiner spaces
$\Hom(\beta_3,\beta_1\beta_2)$ and
$\Hom(\vartheta_+(\beta_3),
\vartheta_+(\beta_1)\vartheta_+(\beta_2))$
are equal.
\endproof

\begin{lemma}
The bijections $\vartheta_\pm$ restrict to
bijections $\MXMo\rightarrow\MXMopm$ of the
ambichiral subsystems.
\labl{ambires}
\end{lemma}

\proof
Let $\tau\in\MXMo$, i.e.\ there are
isometries $s\in\Hom(\tau,\a_\la^+)$ and
$t\in\Hom(\tau,\a_\mu^-)$ for some $\la,\mu\in\NXN$.
Put $q=ts^*\in\Hom(\a^+_\la,\a^-_\mu)$. Then
$q\in\Hom(\iota\la,\iota\mu)$ and
\[ q \eps^+(\la,\rho)^* x_\rho = \eps^-(\mu,\rho)^* x_\rho q \]
whenever $x_\rho\in\Hom(\iota,\iota\rho)$.
Hence, using \erf{anat}, we calculate the left-hand side as
\[ q \eps^+(\la,\rho)^* x_\rho = q \eps^-(\rho,\la) x_\rho
= \eps^- (\rho,\mu) \a^-_\rho(q) x_\rho
= \eps^+(\mu,\rho)^* \a^-_\rho(q) x_\rho \,. \]
Now let us specialize to the case $x_\rho\in\cK_\rho^-$.
Then $x_\rho q=\a^-_\rho(q)x_\rho$, so that
\[ \eps^+(\mu,\rho)^* x_\rho q = \eps^-(\mu,\rho)^* x_\rho q \,. \]
It follows $\alt -+\mu(m)q=\alt --\mu(m)q$ for all $m\in M_-$.
But note that $qq^*=tt^*$ which lies in
$\Hom(\alt --\mu,\alt --\mu)$ by Lemma \ref{Hom=Homt}.
Hence $\alt -+\mu(m)tt^*=\alt --\mu(m)tt^*$  for all $m\in M_-$.
We can similarly derive that
$\alt ++\la(m)ss^*=\alt +-\la(m)ss^*$  for all $m\in M_+$.
There are isometries
$\tilde t\in\Hom(\vartheta_-(\tau),\alt --\mu)$ and
$\tilde s\in\Hom(\vartheta_+(\tau),\alt ++\la)$ such
that $\tilde t \tilde{t}^*=tt^*$ and
$\tilde s \tilde{s}^*=ss^*$. But we now find
\[ \vartheta_-(\tau)(m)=\tilde{t}^* \alt --\mu (m) \tilde t
=\tilde{t}^* \alt -+\mu (m) \tilde t \,,\qquad m\in M_-\,,\]
as well as
\[ \vartheta_+(\tau)(m)=\tilde{s}^* \alt ++\la (m) \tilde s
=\tilde{s}^* \alt +-\mu (m) \tilde s \,,\qquad m\in M_+\,,\]
i.e.\ $\vartheta_\pm(\tau)\in\MXMopm$. Thus
$\vartheta_\pm$ map $\MXMo$ into $\MXMopm$. But it
follows from Proposition \ref{aglob} and Theorem \ref{type1Zs}
that the systems $\MXMo$ and $\MXMopm$ all have the same
ambichiral global index $w_0$. This proves the lemma.
\endproof

We now can state the precise relation between the
coupling matrix $Z$, arising from $N\subset M$ and
given as in \erf{Zb+b-}, and its \typei\ parents
$Z^\pm$ arising from $N\subset M_\pm$.

\begin{theorem}
The entries of the \typei\ coupling matrices $Z^\pm$
arising from $N\subset M_\pm$ can be written as
\begin{equation}
Z^\pm_{\la,\mu}= \sum_{\tau\in\MXMo} b^\pm_{\tau,\la} b^\pm_{\tau,\mu}
\lableq{Z+bb}
with chiral branching coefficients
\begin{equation}
b^\pm_{\tau,\la}= \lan \tau, \a^\pm_\la \ran \,,
\qquad \tau\in\MXMo\,,\quad \la\in\NXN \,.
\lableq{btl}
If the two intermediate subfactors of Corollary \ref{Mpm}
are identical, $M_+=M_-$, (so that the parent coupling
matrices coincide, $Z^+=Z^-$)
then the entries of the coupling matrix $Z$
arising from the full $N\subset M$ can be written as
\begin{equation}
Z_{\la,\mu}= \sum_{\tau\in\MXMo} b^+_{\tau,\la}
b^+_{\omega(\tau),\mu}\,.
\lableq{Zbomb}
Here the permutation $\omega=\vartheta_+^{-1}\circ\vartheta_-$
of $\MXMo$, satisfying $\omega(0)=0$, realizes an automorphism
of the ambichiral fusion rules.
\labl{moseaut}
\end{theorem}

\proof
Since the chiral locality condition holds for $N\subset M_\pm$
we have
\[ \lan\tilde{\tau}_\pm,\alt \pm+\mu\ran
= \lan{\co\iota}_\pm\tilde{\tau}_\pm\iota_\pm,
\mu\ran = \lan\tilde{\tau}_\pm,\alt \pm-\mu\ran \]
for $\tilde{\tau}_\pm\in\MXMopm$ and $\mu\in\NXN$,
thanks to \cite[Prop.\ 3.3]{BE3}. Therefore
\[ \begin{array}{ll}
Z^\pm_{\la,\mu} &= \sum_{\tilde{\tau}_\pm\in\MXMopm}
\lan \alt \pm+\la,\tilde{\tau}_\pm\ran\lan\tilde{\tau}_\pm,
\alt \pm-\mu \ran
=\sum_{\tilde{\tau}_\pm\in\MXMopm} \lan \alt \pm\pm\la,
\tilde{\tau}_\pm\ran
\lan\tilde{\tau}_\pm,\alt \pm\pm\mu \ran \\[.4em]
&=\sum_{\tilde{\tau}_\pm\in\MXMopm}
\lan \a^\pm_\la,\vartheta_\pm^{-1}(\tilde{\tau}_\pm)\ran
\lan\vartheta_\pm^{-1}(\tilde{\tau}_\pm),\a^\pm_\mu \ran
=\sum_{\tau\in\MXMo} b^\pm_{\tau,\la} \DS{b^\pm_{\tau,\mu}}
\end{array} \]
for $\la,\mu\in\NXN$. Now if $M_+=M_-$ then
\[ \begin{array}{ll}
Z_{\la,\mu} &= \sum_{\tau\in\MXMo}
\lan \a^+_\la,\tau\ran\lan\tau,\a^-_\mu \ran
=\sum_{\tau\in\MXMo} b^+_{\tau,\la}
\lan\vartheta_-(\tau),\alt --\mu \ran \\[.4em]
&=\sum_{\tau\in\MXMo} b^+_{\tau,\la}
\lan\vartheta_-(\tau),\alt ++\mu \ran
=\sum_{\tau\in\MXMo} b^+_{\tau,\la}
\lan\vartheta_+^{-1}\circ\vartheta_-(\tau),\a^+_\mu \ran
\end{array} \]
for $\la,\mu\in\NXN$. As $M_+=M_-$, putting
$\om=\vartheta_+^{-1}\circ\vartheta_-$ gives a well-defined
permutation of $\MXMo$ and yields the desired formula
$Z_{\la,\mu}=\sum_{\tau\in\MXMo}b^+_{\tau,\la}b^+_{\om(\tau),\mu}$.
Due to Lemma \ref{furuiso}, $\om$ preserves the fusion rules
and we also have $\om(0)=0$ because always
$\vartheta_\pm(\id_M)=\id_{M_\pm}$.
\endproof

Note that even if $M_+\neq M_-$, $Z^+\neq Z^-$,
the coupling matrix $Z$ is still governed by an
isomorphism of (ambichiral) fusion rule algebras,
in perfect agreement with \cite{MS}.
Namely, if we use the system $\tilde{\tau}_+\in\MXMop$
for the summation in $Z$, then the general formula \erf{Zb+b-}
can be written as
\[ Z_{\la,\mu} = \sum_{\tilde{\tau}_+\in\MXMop}
\lan\alt ++\la,\tilde{\tau}_+\ran
\lan \vartheta (\tilde{\tau}_+),\alt --\mu\ran
= \sum_{\tilde{\tau}_+\in\MXMop}
\lan \la,{\co\iota}_+ \tilde{\tau}_+\iota_+\ran
\lan {\co\iota}_- \vartheta(\tilde{\tau}_+) \iota_-,\mu\ran \]
for $\la,\mu\in\NXN$,  by virtue of Lemma \ref{furuiso} and
chiral locality, guaranteeing \cite[Prop.\ 3.3]{BE3}.
Here $\vartheta$ is the bijection
$\vartheta=\vartheta_-\circ\vartheta_+^{-1}:\MXMop\rightarrow\MXMom$,
yielding an isomorphism of the ambichiral fusion rules.
This corresponds to the ``extended'' coupling matrix
\be
Z^\ext_{\tilde{\tau}_+,\tilde{\tau}_-}=\del {\tilde{\tau}_-}
{\vartheta(\tilde{\tau}_+)} \,, \qquad \tilde{\tau}_+\in\MXMop \,,
\quad \tilde{\tau}_-\in\MXMom \,,
\lableq{Zhetero}
which now has different left and right labels.
(The reader may think of left and right extended characters
$\chi^{\ext;\pm}_{\tilde{\tau}_\pm}=\sum_\la \lan \la,
{\co\iota}_\pm \tilde{\tau}_\pm\iota_\pm\ran \chi_\la$ here.)
This extended coupling matrix also appears in
\erf{bijext} below. Now only if $M_+=M_-$, so that
$\MXMop=\MXMom$, \erf{Zhetero} can be reduced to
$Z^\ext_{\tilde{\tau}_+,\tilde{\tau}_+'}=\del {\tilde{\tau}_+'}
{\vartheta(\tilde{\tau}_+)}$ which is nothing but
the permutation $\om$ up to relabeling by $\MXMo$,
$\om=\vartheta_+^{-1}\circ\vartheta\circ\vartheta_+$,
which we used in \erf{Zbomb} just for notational convenience.

Recall that there is a relative braiding between the
morphisms in $\MXMp$ and $\MXMm$ which restricts to
a proper braiding on $\MXMo$, and for $\tau,\tau'\in\MXMo$
these braiding operators are given by \cite{BE3}
\[ \epsr(\tau,\tau') = s^* \a^-_\mu(t)^* \eps^+(\la,\mu)
\a^+_\la (s)t \]
whenever $t\in\Hom(\tau,\a^+_\la)$ and
$s\in\Hom(\tau',\a^-_\mu)$ are isometries,
$\la,\mu\in\NXN$. We can extend this braiding from $\MXMo$
to $\Sigma(\MXMo)$ as explained in \cite[Sect.\ 2]{BEK1}.
Now let $N^\op$ denote the opposite algebra of $N$ and
let $j$ denote the natural anti-linear isomorphism. For
$\la\in\End(N)$ we denote $\la^\op=j\circ\la\circ j$.
We proceed analogously for $M_-^\op$, the opposite algebra
of $M_-$.

\begin{proposition}
There exists a (\typeiii) factor $B$ such that we have
irreducible inclusions
\be
N\otimes N^\op \subset M_+ \otimes M_-^\op \subset B
\ee
with the following properties:
\begin{enumerate}
\item The dual canonical endomorphism $\Theta_\ext$ of the
 inclusion $M_+\otimes M_-^\op \subset B$ decomposes as
\be
[\Theta_\ext] =\bigoplus_{\tau\in\MXMo}
[\vartheta_+(\tau) \otimes \vartheta_-(\tau)^\op] \,.
\lableq{bijext}
\item The dual canonical endomorphism $\Theta$ of the
 inclusion $N\otimes N^\op \subset B$ decomposes as
\be
[\Theta] =\bigoplus_{\la,\mu\in\NXN}
Z_{\la,\mu} \, [\la \otimes \mu^\op] \,.
\lableq{Zcft}
\item If $\tau_\pm,\sig_\pm\in\End(M)$ are the extensions
 of $\tilde{\tau}_\pm,\tilde{\sig}_\pm\in\MXMopm$,
 respectively, according to Lemma \ref{bij}, then
\be
y \sig_+(x) \epsr (\tilde{\tau}_+,\tilde{\sig}_+)
= \epsr (\tilde{\tau}_-,\tilde{\sig}_-) \tau_-(y)x
\lableq{extnat}
holds whenever $x\in\Hom(\tau_+,\tau_-)$
and $y\in\Hom(\sig_+,\sig_-)$.
\end{enumerate}
\labl{extcft}
\end{proposition}

\proof
Lemma \ref{bij} (together with Lemma \ref{ambires}) tells us
that $\tilde{\tau}_\pm\in\MXMopm$ can be extended to
$\tau_\pm\in\End(M)$ such that $\tau_+$ and $\tau_-$ are
equivalent to some morphisms in $\MXMo$, and we have
$[\tau_+]=[\tau_-]$ if and only if
$\tilde{\tau}_\pm=\vartheta_\pm(\tau)$ for a $\tau\in\MXMo$.
Using these extensions for subfactors $M_\pm\subset M$
with systems $\MXMopm$, then \cite{R7} determines a factor
$B$ such that $M_+\otimes M_-^\op\subset B$ is an irreducible
subfactor with its dual canonical endomorphism $\Theta_\ext$
decomposing as
\[ [\Theta_\ext] = \bigoplus_{\tilde{\tau}_+\in\MXMop}
\bigoplus_{\tilde{\tau}_-\in\MXMom} \lan\tau_+,\tau_-\ran
[\tilde{\tau}_+ \otimes \tilde{\tau}_-^\op]=
\bigoplus_{\tau\in\MXMo}
[\vartheta_+(\tau) \otimes \vartheta_-(\tau)^\op] \,, \]
proving {\it 1}. Now note that the injection map for
$N\otimes N^\op \hookrightarrow M_+\otimes M_-^\op$
is given by $\iota_+\otimes\iota_-^\op$. Therefore
the dual canonical endomorphism for
$N\otimes N^\op\subset B$ is obtained as
$\Theta=({\co\iota}_+\otimes{\co\iota}_-^\op)
\circ\Theta^\ext\circ(\iota_+\otimes\iota_-^\op)$
so that
\[ [\Theta] = \bigoplus_{\tau\in\MXMo}
[({\co\iota}_+\circ\vartheta_+(\tau)\circ\iota_+)\otimes
({\co\iota}_-\circ\vartheta_-(\tau)\circ\iota_-)^\op]\,.\]
Now
\[ [{\co\iota}_\pm\circ\vartheta_\pm(\tau)\circ\iota_\pm]
=\bigoplus_{\la\in\NXN} \lan {\co\iota}_\pm\circ\vartheta_\pm
(\tau)\circ\iota_\pm,\la\ran [\la] \,, \]
and since the subfactors $N\subset M_\pm$ satisfy
chiral locality one has
\[ \lan {\co\iota}_\pm\circ\vartheta_\pm (\tau)
\circ\iota_\pm,\la\ran = \lan \vartheta_\pm(\tau),
\alt \pm\pm \la \ran = \lan \tau, \a^\pm_\la \ran
= b_{\tau,\la}^\pm \]
by virtue of \cite[Prop.\ 3.3]{BE3} and Lemma \ref{furuiso}.
Hence
\[ [\Theta] = \bigoplus_{\la,\mu\in\NXN}
\bigoplus_{\tau\in\MXMo} b_{\tau,\la}^+ b_{\tau,\mu}^-
[\la \otimes \mu^\op] = \bigoplus_{\la,\mu\in\NXN}
Z_{\la,\mu} [\la \otimes \mu^\op] \,, \]
proving {\it 2}.
Finally, if $\tau_\pm,\sig_\pm\in\End(M)$ denote the
extensions of $\tilde{\tau}_\pm,\tilde{\sig}_\pm\in\MXMopm$,
respectively, as in Lemma \ref{bij}, then there are some
$\la_+,\la_-,\mu_+,\mu_-\in\NXN$ and isometries
$\tilde{t}_\pm\in\Hom(\tilde{\tau}_\pm,\alt \pm\pm{\la_\pm})$
and 
$\tilde{s}_\pm\in\Hom(\tilde{\sig}_\pm,\alt \pm\pm{\mu_\pm})$
so that
$\tau_\pm(m)=\tilde{t}_\pm^* \a^\pm_{\la_\pm}(m) \tilde{t}_\pm$
and 
$\sig_\pm(m)=\tilde{s}_\pm^* \a^\pm_{\mu_\pm}(m) \tilde{s}_\pm$
for all $m\in M$. Note that also
$\tilde{t}_\pm\in\Hom(\tilde{\tau}_\pm,\alt \pm\mp{\la_\pm})$
and 
$\tilde{s}_\pm\in\Hom(\tilde{\sig}_\pm,\alt \pm\mp{\mu_\pm})$
because chiral locality holds for $N\subset M_\pm$ and then
ambichiral morphisms are obtained from
$\a^+$- and $\a^-$-induction by use of the same isometries
\cite[Sect.\ 3]{BE3}.
Hence we have
\[ \begin{array}{ll}
\epsr (\tilde{\tau}_\pm,\tilde{\sig}_\pm) &=
\tilde{s}_\pm^* \alt \pm\pm{\mu_\pm}(\tilde{t}_\pm)^*
\eps^+(\la_\pm,\mu_\pm)\alt \pm\pm{\la_\pm}(\tilde{s}_\pm)
\tilde{t}_\pm \\[.4em]
&= \tilde{s}_\pm^* \a^\pm_{\mu_\pm}(\tilde{t}_\pm)^*
\eps^+(\la_\pm,\mu_\pm)\a^\pm_{\la_\pm}(\tilde{s}_\pm)
\tilde{t}_\pm \,,
\end{array} \]
where we also used Corollary \ref{res}. We now can compute
\[ \begin{array}{ll}
y\sig_+(x) \epsr (\tilde{\tau}_+,\tilde{\sig}_+)
&= y\sig_+(x) \tilde{s}_+^* \a^+_{\mu_+}(\tilde{t}_+)^*
\eps^+(\la_+,\mu_+)\a^+_{\la_+}(\tilde{s}_+) \tilde{t}_+ \\[.4em]
&= y\sig_+(x \tilde{t}_+^*) \tilde{s}_+^* \eps^+(\la_+,\mu_+)
\a^+_{\la_+}(\tilde{s}_+) \tilde{t}_+ \\[.4em]
&= \sig_-(x \tilde{t}_+^*) \tilde{s}_-^*\tilde{s}_- y \tilde{s}_+^*
\eps^+(\la_+,\mu_+) \a^+_{\la_+}(\tilde{s}_+) \tilde{t}_+ \\[.4em]
&= \tilde{s}_-^* \a^-_{\mu_-}
(\tilde{t}_-^*\tilde{t}_- x \tilde{t}_+^*)
\eps^+(\la_+,\mu_-) \a^+_{\la_+}(\tilde{s}_- y)\tilde{t}_+ \\[.4em]
&= \tilde{s}_-^* \a^-_{\mu_-} (\tilde{t}_-^*)
\eps^+(\la_-,\mu_-) \tilde{t}_- x \tilde{t}_+^*
\a^+_{\la_+}(\tilde{s}_-)\tilde{t}_+ \tau_+(y) \\[.4em]
&= \tilde{s}_-^* \a^-_{\mu_-} (\tilde{t}_-^*)
\eps^+(\la_-,\mu_-)  \a^-_{\la_-}(\tilde{s}_-)
\tilde{t}_- x \tau_+(y)
= \epsr (\tilde{\tau}_-,\tilde{\sig}_-) \tau_-(y)x \,,
\end{array} \]
where we used \erf{anat} twice, proving {\it 3}.
\endproof

The relevance of Proposition \ref{extcft} is the following.
Suppose that our factor $N$ is obtained as a local factor
$N=N(I_\circ)$ of a quantum field theoretical net of factors
$\{N(I)\}$ indexed by proper intervals $I\subset \bbR$
on the real line, and that the system $\NXN$ is obtained
as restrictions of DHR-morphisms (cf.\ \cite{H}) to $N$.
This is in fact the case in our RCFT examples arising from
current algebras where the net is defined in terms of local
loop groups in the vacuum representation.
Taking two copies of such a net and placing the real axes
on the light cone, then this defines a local net $\{A(\cO)\}$,
indexed by double cones $\cO$ on two-dimensional
Minkowski space (cf.\ \cite{R5} for such constructions).
Given a subfactor $N\subset M$, determining in turn two
subfactors $N\subset M_\pm$ obeying chiral locality,
will provide two local nets of subfactors
$\{N(I)\subset M_\pm(I) \}$ due to \cite{LR}.
Arranging $M_+(I)$ and $M_-(J)$ on the two light
cone axes defines a local net of subfactors
$\{A(\cO)\subset A_\ext(\cO)\}$ in Minkowski space.
The embedding $M_+\otimes M_-^\op \subset B$
gives rise to another net of subfactors
$\{A_\ext(\cO) \subset B(\cO)\}$, and \erf{extnat}
ensures that the net $\{B(\cO)\}$ satisfies locality,
due to Rehren's recent result \cite{R7}.
As already shown in \cite{R7}, there exist a local
two-dimensional quantum field theory such that the
coupling matrix $Z$ describes its restriction to
the tensor products of its chiral building blocks $N(I)$,
and this is here expressed in \erf{Zcft}.
Now \erf{bijext} tells us that there are chiral extensions
$N(I)\subset M_+(I)$ and $N(I)\subset M_-(I)$ for left and
right chiral nets which are indeed maximal as the coupling
matrix for $\{A_\ext(\cO) \subset B(\cO)\}$ is a bijection.
This shows that the inclusions $N\subset M_\pm$ should in
fact be regarded as the subfactor version of left- and
right maximal extensions of the chiral algebra.

\section{Extended S- and T-matrices}
\labl{SText}

Using the braiding arising from the relative braiding one can
define the statistics phase $\om_\tau$ of $\tau\in\MXMo$ by
$d_\tau\phi_\tau (\epsr(\tau,\tau))=\om_\tau\bfe$.

\begin{lemma}
Let $\tau\in\MXMo$ such that $[\tau]$ is a subsector of
$[\a^+_\la]$ and $[\a^-_\mu]$ for some $\la,\mu\in\NXN$.
Then we have $\om_\la=\om_\tau=\om_\mu$.
\labl{omomom}
\end{lemma}

\proof
Let $t\in\Hom(\tau,\a^+_\la)$ and $s\in\Hom(\tau,\a^-_\mu)$
be isometries. Then
\[ \begin{array}{ll}
\om_\tau \bfe &= d_\tau \phi_\tau (\epsr(\tau,\tau))
= d_\tau \phi_\tau (s^* \a^-_\mu(t)^* \eps^+(\la,\mu)
\a^+_\la (s)t) \\[.4em]
&= d_\tau \phi_\tau (\tau(t)^* s^* \eps^+(\la,\mu) t \tau(s))
= d_\tau t^* \phi_\tau (s^* \eps^+(\la,\mu) t)s \\[.4em]
&= d_\la t^* \phi_{\a^+_\la} (ts^* \eps^+(\la,\mu))s
= d_\la t^* \phi_{\a^+_\la}
(\eps^+(\la,\la)\a^+_\la(ts^*))s  \\[.4em]
&= d_\la t^* \phi_{\a^+_\la} (\eps^+(\la,\la))ts^*s \,,
\end{array} \]
where we used \erf{phiaba}, \erf{phibttb}, and since
$ts^*\in\Hom(\iota\mu,\iota\la)$ we could also apply
\erf{anat}. Note that $\phi_{\a^+_\la}$ can be given
as $\phi_{\a^+_\la}(m)=r_\la^*\a^+_{\co\la}(m)r_\la$
for all $m\in M$, so that in particular
$\phi_{\a^+_\la}(n)=\phi_\la(n)$ for $n\in N$.
Hence $\om_\tau=\om_\la$. We can compute similarly
\[ d_\tau t^* \phi_\tau (s^* \eps^+(\la,\mu) t)s =
d_\mu t^* \phi_{\a^-_\mu} (\eps^+(\la,\mu) ts^*)s =
d_\mu t^* \phi_{\a^-_\mu}
(\a^-_\mu(ts^*)\eps^+(\mu,\mu))s \,,\]
establishing $\om_\tau=\om_\mu$.
\endproof

Note that with the expansion \erf{Zb+b-},
Lemma \ref{omomom} implies easily
$\om_\la Z_{\la,\mu}=Z_{\la,\mu} \om_\mu$, i.e.\ it
gives a new and simple proof of the commutativity of 
the matrices $\Omega$ and $Z$ which was first
established in \cite[Thm.\ 5.7]{BEK1}.

\begin{lemma}
For $\beta\in\MXM$ we have
\be
\sum_{\la\in\NXN} \om_\la d_\la \lan \beta,\a^\pm_\la \ran
= \left\{ \begin{array}{c@{\qquad:\qquad}l}
 \om_\tau d_\tau w/w_+ & \beta=\tau\in\MXMo \\[.6em]
 0 & \mbox{otherwise.} \end{array} \right.
\label{inphsum}
\ee
\labl{inphase}
\end{lemma}

\proof
All we need to show is that the left hand side of \erf{inphsum}
vanishes whenever $\beta\notin\MXMo$ because we recall once
more from (the proof of) \cite[Prop.\ 3.1]{BEK2} that we have
$d_\tau w/w_+ =\sum_\la d_\la b^\pm_{\tau,\la}$
for any $\tau\in\MXMo$, and then the claim follows from
Lemma \ref{omomom}. For this purpose we define
vectors $\vec{u}^\pm$ with entries
\[ u^\pm_\beta=\sum_{\la\in\NXN} \om_\la d_\la
\lan\beta,\a^\pm_\la\ran \,, \qquad \beta\in\MXM \,.\]
We clearly have
\[ \begin{array}{ll}
\| \vec{u}^+ \|^2 &= \sum_\beta \sum_{\la,\nu} \om_\la \om_\nu^{-1}
 d_\la d_\nu \lan \a^+_\nu,\beta \ran \lan \beta, \a^+_\la \ran
= \sum_{\la,\nu} \om_\la \om_\nu^{-1}
 d_\la d_\nu \lan \a^+_{\co\la} \a^+_\nu, \id \ran \\[.4em]
&= \sum_{\la,\mu,\nu} \om_\la \om_\nu^{-1}
 d_\la d_\nu N_{\co\la,\nu}^\mu \lan \a^+_\mu, \id \ran
= \sum_{\la,\mu,\nu} \om_\la \om_\nu^{-1}
 d_\la d_\nu N_{\la,\mu}^\nu \om_\mu Z_{\mu,0} \\[.4em]
&= \sum_{\la,\mu} Y_{0,\la} Y_{\la,\mu} Z_{\mu,0} \,,
\end{array} \]
where we used that $Z_{\mu,0}=\om_\mu Z_{\mu,0}$
by Lemma \ref{omomom}. Similarly we obtain
$\|\vec{u}^-\|^2=\sum_{\la,\mu}Y_{0,\la}Y_{\la,\mu}Z_{0,\mu}$.
On the other hand, we can compute the inner product
\[ \begin{array}{ll}
\lan \vec{u}^+,\vec{u}^- \ran
&= \sum_{\la,\mu} \sum_\beta \om_\la^{-1} d_\la d_\mu \om_\mu
\lan \a^+_\la,\beta\ran \lan\beta,\a^-_\mu \ran
= \sum_{\la,\mu} \om_\la^{-1} d_\la d_\mu \om_\mu
 Z_{\la,\mu} \\[.4em]
&= \sum_{\la,\mu} d_\la Z_{\la,\mu} d_\mu
= \sum_{\la,\mu} Y_{0,\la} Z_{\la,\mu} Y_{\mu,0}
= \sum_{\la,\mu} Y_{0,\la} Y_{\la,\mu} Z_{\mu,0}
= \| \vec{u}^+ \|^2\,,
\end{array} \]
where we used the commutativity of $\Omega$ and $Y$
with $Z$ of \cite[Thm.\ 5.7]{BEK1}.
Commuting $Z$ in the fifth equality to the left rather
than to the right gives
$\lan\vec{u}^+,\vec{u}^-\ran=\|\vec{u}^-\|^2$
Thus we conclude $\vec{u}^+=\vec{u}^-$. Since obviously
$u^\pm_\beta=0$ whenever $\beta\notin\MXMpm$ this
implies $u^\pm_\beta=0$ whenever $\beta\notin\MXMo$.
\endproof

Let $Y^\ext$ and $\Omega^\ext$ denote the Y- and
$\Omega$-matrices associated to the braided system $\MXMo$.

\begin{lemma}
We have
\be
\frac w{w_+} \sum_{\tau'\in\MXMo}
Y^\ext_{\tau,\tau'} b^\pm_{\tau',\mu} =
\sum_{\la\in\NXN} b^\pm_{\tau,\la} Y_{\la,\mu} \,,
\qquad \sum_{\tau'\in\MXMo}
\Omega^\ext_{\tau,\tau'} b^\pm_{\tau',\mu}
= \sum_{\la\in\NXN} b^\pm_{\tau,\la} \Omega_{\la,\mu}
\lableq{Yextbb}
for all $\tau\in\MXMo$ and $\mu\in\NXN$.
\labl{YOmextbYOm}
\end{lemma}

\proof
Note that the second relation in \erf{Yextbb} is nothing
but Lemma \ref{omomom}, as this is just
$\om_\tau b^\pm_{\tau,\la} = b^\pm_{\tau,\la} \om_\la$.
So we just need to verify the first relation in \erf{Yextbb}.
We now compute
\[ \begin{array}{ll}
\DS{\sum_\la} b^\pm_{\tau,\la} Y_{\la,\mu}
&= \DS{\sum_{\la,\rho} \frac{\om_\la \om_\mu}{\om_\rho}}
N_{\la,\mu}^\rho d_\rho b^\pm_{\tau,\la} 
= \DS{\sum_{\la,\rho} \frac{\om_\tau \om_\mu}{\om_\rho}}
N_{\rho,\co\mu}^\la d_\rho \lan\a^\pm_\la, \tau \ran \\[1.4em]
&= \DS{\sum_\rho \frac{\om_\tau \om_\mu}{\om_\rho}} d_\rho
\lan\a^\pm_\rho, \tau \a^\pm_\mu\ran 
= \DS{\sum_\rho \sum_\beta \frac{\om_\tau \om_\mu}{\om_\rho}}
 d_\rho \lan\a^\pm_\rho, \beta \ran \lan \beta,
 \tau \a^\pm_\mu\ran \\[.4em]
&= \DS{\frac w{w_+} \sum_{\tau''}
 \frac{\om_\tau \om_\mu}{\om_{\tau''}}}
  d_{\tau''} \lan \tau'',\tau \a^\pm_\mu\ran
= \DS{\frac w{w_+} \sum_{\tau',\tau''}
 \frac{\om_\tau \om_\mu}{\om_{\tau''}}} N_{\co\tau,\tau''}^{\tau'}
  d_{\tau''} \lan \tau', \a^\pm_\mu\ran \\[.4em]
&= \DS{\frac w{w_+} \sum_{\tau',\tau''}
 \frac{\om_\tau \om_{\tau'}}{\om_{\tau''}}}
 N_{\tau,\tau'}^{\tau''} d_{\tau''} b^\pm_{\tau',\mu}
= \DS{\frac w{w_+} \sum_{\tau'}} Y^\ext_{\tau,\tau'} 
 b^\pm_{\tau',\mu} \,,
\end{array} \]
where we used (the complex conjugate of)
Lemma \ref{inphase} in the fifth equality.
\endproof

Recall from Section \ref{prelim} the complex number
$z=\sum_{\la\in\NXN} d_\la^2 \om_\la$.
Analogously we put
$z_0=\sum_{\tau\in\MXMo} d_\tau^2 \om_\tau$.

\begin{lemma}
We have $z_0= \DS{\frac {w_+}w} z$.
\labl{c=cext}
\end{lemma}

\proof
Using Lemma \ref{inphase} we can compute
\[ z_0 = \sum_{\tau\in\MXMo} \om_\tau d_\tau^2 
= \frac {w_+}w \sum_{\beta\in\MXM} \sum_{\la\in\NXN}
 \om_\la d_\la \lan \beta,\a^\pm_\la \ran d_\beta
= \frac {w_+}w \sum_{\la\in\NXN} \om_\la d_\la^2 
= \frac {w_+}w z \,, \]
where we used
$\sum_{\beta\in\MXM}\lan\beta,\a^\pm_\la\ran d_\beta=d_\la$.
\endproof

Assume that $z\neq 0$ so that the central charge can be
defined. Since $w/w_+$ is a real number, Lemma \ref{c=cext}
tells us that the central charges $c$ and $c_\ext$
of the braided systems $\NXN$ and $\MXMo$, respectively,
which are defined modulo 8, coincide.
As a corollary of Lemma \ref{YOmextbYOm}, the intertwining
properties of chiral branching coefficients between original
and extended S- and T-matrices are therefore clarified in
the following

\begin{theorem}
Provided $z\neq 0$ one has
\be
\sum_{\tau'\in\MXMo}
S^\ext_{\tau,\tau'} b^\pm_{\tau',\mu} =
\sum_{\la\in\NXN} b^\pm_{\tau,\la} S_{\la,\mu} \,,
\qquad \sum_{\tau'\in\MXMo}
T^\ext_{\tau,\tau'} b^\pm_{\tau',\mu}
= \sum_{\la\in\NXN} b^\pm_{\tau,\la} T_{\la,\mu}
\lableq{Sextbb}
for all $\tau\in\MXMo$ and $\mu\in\NXN$.
\labl{STextbST}
\end{theorem}

We would like to remind the reader that,
if the braiding on $\NXN$ is non-degenerate,
then so is the ambichiral braiding \cite[Thm.\ 4.2]{BEK2}.
In other words, whenever the original S- and T-matrices
are modular then so are the extended S- and T-matrices.

Now let us return to our situation $N\subset M_\pm \subset M$
and apply the above results also to the subfactors
$N\subset M_\pm$. Let $Y^{\ext;\pm}$
and $\Omega^{\ext;\pm}$ the Y- and $\Omega$-matrices
associated to the braided systems $\MXMopm$.
Recalling now $Z^+_{\la,0}=Z_{\la,0}$ and
$Z^-_{0,\la}=Z_{0,\la}$,
we obtain from Lemmata \ref{furuiso}, \ref{ambires},
\ref{omomom} and \ref{c=cext} the following

\begin{theorem}
The matrices $Y^\ext$, $\Omega^\ext$,
and $Y^{\ext;\pm}$, $\Omega^{\ext;\pm}$
coincide subject to the bijections $\vartheta_\pm$.
If $z\neq 0$, then so do the corresponding S- and T-matrices
which are then well-defined. In formulae,
\be
S^\ext_{\tau,\tau'}
= S^{\ext;\pm}_{\vartheta_\pm(\tau),\vartheta_\pm(\tau')}\,,
\qquad  T^\ext_{\tau,\tau'}
= T^{\ext;\pm}_{\vartheta_\pm(\tau),\vartheta_\pm(\tau')}\,,
\ee
for all $\tau,\tau'\in\MXMo$.
\labl{SSTT}
\end{theorem}

We remark that in the case that $M_+=M_-$ one obtains
by use of the properties of the relative braiding operators
\cite[Lemma 3.11]{BE3} and from Corollary \ref{res},
that the ambichiral braiding operators are the same for
$\MXMo$ and $\MXMopm$, subject to the bijections
$\vartheta_\pm$, so that Theorem \ref{SSTT} is trivial
in this case.

\section{Heterotic examples}
\labl{heterotic}

We consider the $\SOn$ current algebra models at level 1, and
where $n$ is a multiple of 16, $n=16\ell$, $\ell=1,2,3,...\,$.
These theories have four sectors, the basic ($0$), vector (v),
spinor (s) and conjugate spinor (c) module, corresponding to
highest weights $0$, $\Lambda_{(1)}$, $\Lambda_{(r-1)}$
and $\Lambda_{(r)}$, respectively; here
$r=n/2=8\ell$ is the rank of $\SOn$. The conformal
dimensions are given as $h_0=0$, $h_\rmv=1/2$,
$h_\rms=h_\rmc=\ell$, and the sectors obey
$\bbZ_2\times\bbZ_2$ fusion rules. The
Kac-Peterson matrices are given explicitly as
\be
S = \frac 12 \left( \begin{array}{rrrr}
1 & 1 & 1 & 1 \\ 1 & 1 & -1 & -1 \\
1 & -1 & 1 & -1 \\ 1 & -1 & -1 & 1
\end{array} \right) , \qquad
T = \E^{-2\pi\I \ell/3} \left( \begin{array}{rrrr}
1 & 0 & 0 & 0 \\ 0 & -1 & 0 & 0 \\
0 & 0 & 1 & 0 \\ 0 & 0 & 0 & 1
\end{array}  \right) .
\lableq{ST}
It is easy to check that there are exactly six
modular invariants, $Z=\bfe$, $W$, $X_\rms$,
$X_\rmc$, $Q$, $\tmat Q$. Here
\[ W = \left( \begin{array}{rrrr}
1 & 0 & 0 & 0 \\ 0 & 1 & 0 & 0 \\
0 & 0 & 0 & 1 \\ 0 & 0 & 1 & 0
\end{array}  \right) \,, \qquad
X_\rms = \left( \begin{array}{rrrr}
1 & 0 & 1 & 0 \\ 0 & 0 & 0 & 0 \\
1 & 0 & 1 & 0 \\ 0 & 0 & 0 & 0
\end{array}  \right) \,, \qquad
Q = \left( \begin{array}{rrrr}
1 & 0 & 0 & 1 \\ 0 & 0 & 0 & 0 \\
1 & 0 & 0 & 1 \\ 0 & 0 & 0 & 0
\end{array}  \right) \,, \]
and $X_\rmc=W X_\rms W$.
(Note that $Q=X_\rms W$ and $\tmat Q=W X_\rms$.)
The matrix $Q$ and its transpose $\tmat Q$
are two examples for modular invariants with
non-symmetric vacuum coupling. Such ``heterotic''
invariants seem to be extremely rare and have not
enjoyed particular attention in the literature,
perhaps because they were erroneously dismissed
as being spurious in the sense that they would not
correspond to a physical partition function.
Examples for truly spurious modular invariants
were given in \cite{SY,V,FSS} and found
to be ``coincidental'' linear combinations of
proper physical invariants. Note that although
there is a linear dependence here, namely
\[ \bfe - W - X_\rms - X_\rmc + Q + \tmat Q = 0 \,, \]
we cannot express $Q$ (or $\tmat Q$) alone as a linear
combination of the four symmetric invariants.
This may serve as a first indication that $Q$ and
$\tmat Q$ are not spurious. We will now demonstrate
that they can be realized from subfactors.

The $\bbZ_2\times\bbZ_2$ fusion rules for these models
were proven in the DHR framework in \cite{B2}, and
together with the conformal spin and statistics theorem
\cite{FG1,FRS2,GL} we conclude that there is a net of
\typeiii\ factors on $S^1$ with a system
$\{\id,\rho_\rmv,\rho_\rms,\rho_\rmc\}$ of localized
and transportable, hence braided endomorphisms, such that
the statistics S- and T-matrices are given by \erf{ST}.
Because the statistics phases are given as
$\omega_\rmv=-1$ and $\omega_\rms=\omega_\rmc=1$,
we can assume that the morphisms in the system
obey the $\bbZ_2\times\bbZ_2$ fusion rules even
by individual multiplication,
\[ \rho_\rmv^2=\rho_\rms^2=\rho_\rmc^2=\id\,, \qquad
\rho_\rmv \rho_\rms = \rho_\rms \rho_\rmv = \rho_\rmc \,,\]
thanks to \cite[Lemma 4.4]{R3}.
This is enough to proceed with the DHR construction of
the field net \cite{DHR2}, as already carried out
similarly for simple current extensions with cyclic
groups in \cite{BE2,BE3}. In fact, all we need to do here
is to pick a single local factor $N=N(I)$ such that the
interval $I\subset S^1$ contains the localization region
of the morphisms, and then we construct the cross product
subfactor $N\subset N\rtimes(\bbZ_2\times\bbZ_2)$.
Then the corresponding dual canonical endomorphism
$\canr$ decomposes as a sector as
\[ [\canr] = [\id] \oplus [\rho_\rmv] \oplus
[\rho_\rms] \oplus [\rho_\rmc] \,. \]
Checking
$\lan\iota\la,\iota\mu\ran=\lan\canr\la,\mu\ran=1$
for $\la,\mu=\id,\rho_\rmv,\rho_\rms,\rho_\rmc$,
we find that there is only a single $M$-$N$ sector,
namely $[\iota]$. By \cite[Cor.\ 6.13]{BEK1} we
conclude that the modular invariant coupling matrix $Z$
arising from this subfactor must fulfill $\tr Z=1$.
This leaves only the possibility that $Z$ is $Q$
or $\tmat Q$. We may and do assume that $Z=Q$,
otherwise we exchange braiding and opposite braiding.
It is easy to determine the intermediate subfactors
$N\subset M_\pm \subset M$. Namely, we have
$M_+=N\rtimes_{\rho_\rms} \bbZ_2$ and
$M_-=N\rtimes_{\rho_\rmc} \bbZ_2$ with dual
canonical endomorphism sectors
$[\canr_+]=[\id]\oplus[\rho_\rms]$ and
$[\canr_-]=[\id]\oplus[\rho_\rmc]$, respectively.
That both extensions are local can also be
checked from $\omega_\rms=\omega_\rmc=1$. We therefore
find $Z^+=X_\rms$ and $Z^-=X_\rmc$. Finally, the
permutation invariant $W$ is obtained from the
non-local extension $M_\rmv=N\rtimes_{\rho_\rmv}\bbZ_2$.

\section{Conclusions}

We studied the structure of coupling matrices $Z$
arising from braided subfactors $N\subset M$ through
intermediate subfactors $N\subset M_\pm\subset M$ which in
turn determine \typei\ ``parent'' coupling matrices $Z^\pm$.
We demonstrated that the inclusions $N\subset M_+$ and
$N\subset M_-$ should be recognized as the subfactor
version of left and right maximal chiral algebra extensions.
The main application we have in mind is RCFT where the
S- and T-matrices arising from the braiding are modular.
For current algebra models based on Lie groups $\SUn$ or
others, the coupling matrices from subfactors are then
modular invariants for their Kac-Peterson matrices.
Most but not all of the known modular invariant coupling
matrices of such models are either \typei, \erf{ty1},
or \typeii, \erf{ty2}, and the \typeii\ invariants have a
unique \typei\ parent. For example, the parents of the $\SUz$
\typeii\ modular invariants D$_{2\ell+1}$ ($\ell=2,3,\ldots$)
and E$_7$ are A$_{4\ell-1}$ and D$_{10}$, respectively.
For such invariants the extended left and right chiral
algebras are the same.
(In the D$_{\mathrm{odd}}$ examples the extended algebras
are the original, identical left and right current algebras.)
In fact, the E$_7$ modular invariant has been
constructed from a subfactor with
$[\canr]=[\la_0]\oplus[\la_8]\oplus[\la_{16}]$
in \cite{BEK2}, and here we obtain $M_+=M_-$ with
$[\canr_\pm]=[\la_0]\oplus[\la_{16}]$ which produces
the simple current extension D$_{10}$ invariant.
For the cases D$_{\mathrm{even}}$, E$_6$ and E$_8$
treated in \cite{BE2,BE3} where $N\subset M$ is
subject to chiral locality from the beginning,
we clearly find $M=M_+=M_-$ which indeed are the
local factors of the local chiral extensions
considered e.g.\ in \cite{RST}.
(In fact all invariants obtained from subfactors
obeying chiral locality are clearly their own parents
due to Proposition \ref{condloc}.)
We showed that $Z$ is in fact \typeii, \erf{ty2},
whenever the extensions coincide, $M_+=M_-$.

It is interesting that all our results could be
derived without assuming the non-degeneracy of the
braiding, i.e.\ all our statements are true even if
the modular group is not around. We similarly derived
in \cite{BEK2} without such condition that trivial
vacuum coupling, $Z_{\la,0}=\del \la0$, is equivalent
to $Z$ being a fusion rule automorphism
(and to $Z_{0,\la}=\del \la0$), thus recovering a result
previously encountered in RCFT \cite{DV,MS}.
In this paper we started with a braided subfactor producing
some coupling matrix with possibly non-trivial vacuum
coupling, and our results show that the ``extended''
coupling matrix, \erf{Zhetero}, is a bijection
$\MXMop\rightarrow\MXMom$ which yields an isomorphism 
of the fusion rules of the ambichiral systems.
Moreover, the corresponding extended S- and T-matrices
coincide subject to this isomorphism (Theorem \ref{SSTT}).
In the (modular) RCFT case, they are recognized as the
S- and T-matrices of the extended left and right
chiral algebra, and therefore Theorem \ref{SSTT} provides
in particular a subfactor version of \cite[Eq.\ (4.5)]{MS}.
But note that the derivations of the fusion rule automorphism
in \cite{DV,MS} in turn rely on the Verlinde formula \cite{Ve}
whereas our derivation holds even if the braiding is
degenerate, i.e.\ even if the Verlinde formula does not hold.
Our result comes in the same spirit as \cite{R5} where
the embedding of left and right chiral observables in
a 2D conformal quantum field theory is analyzed and the
corresponding coupling matrix is shown to describe an
automorphism of fusion rules if and only if the
chiral observables are maximal.
Namely, the result of \cite{R5} is derived under very general
assumptions in the framework of local quantum physics,
and it is in particular entirely independent of
the $\SLZ$ machinery heavily exploited in \cite{DV,MS}.

Note that ``almost all'' known modular invariants satisfy
$Z_{\la,0}=Z_{0,\la}$. This means $[\canr_+]=[\canr_-]$,
and this comes close to $M_+=M_-$. In particular,
if $\lan\canr,\la\ran\le 1$, then this forces
$\cK_\la^+=\cK_\la^-$ so that we necessarily have $M_+=M_-$.
Similarly a total degenerate braiding gives
rise to $M_+=M_-$. Nevertheless our example in
Section \ref{heterotic} has shown that $M_+\neq M_-$ is
possible, an that even $Z^+\neq Z^-$ can occur.
The significance of different left and right chiral
extensions reflected in the possibility $M_+\neq M_-$
and even in different parent coupling matrices may
come a little surprising. For example, the related
``heterotic'' extensions of current algebra models are
not particularly well studied objects. One reason may
be that the most popular models, those based on $\SUn$
current algebras, only seem to have modular invariants
with identical parents --- it is in fact likely that all
$\SUn$ invariants are entirely symmetric.\footnote{This
was pointed out to us by Terry Gannon.} But can it happen
that $Z^+=Z^-$ but $M_+\neq M_-$? We do not know an
example but we neither see a reason why this should
not be possible. For instance, if there is
$\lan\la,\canr\ran\ge2$ for some $\la$ then it may happen
that $\cK^+_\la\neq\cK^-_\la$ though these spaces may still
have the same dimension, $Z_{\la,0}=Z_{0,\la}$.
In other words, it is conceivable that certain
modular invariants look like being \typei\ or \typeii\
though they really come from heterotic extensions.

Let us finally mention that the exotic modular invariants
which are argued not to correspond to any RCFT
in \cite{SY,V,FSS}, will not be produced from subfactors
by the machinery of \cite{BEK1,BEK2}.
Note that the standard argument showing that a modular
invariant $Z$ does not give a partition function of a RCFT
is to disprove the existence of an extended S-matrix.
However, from braided subfactors there always arises a
matrix $S^\ext$ with all the required properties.
And in fact, Rehren's recent result \cite{R7} (and in turn
our Proposition \ref{extcft}) shows generally that all
coupling matrices which arise from an embedding of some
local algebra of a chiral RCFT describe the restriction
of a 2D RCFT to its chiral building blocks.
This implies for example that the heterotic modular invariants
$Q$ and $\tmat Q$ discussed in Section \ref{heterotic}
are in fact coupling matrices of a RCFT whose chiral algebras
are different maximal extensions of the $\SOn_1$ current
algebra.

\vspace{0.2cm}\addtolength{\baselineskip}{-2pt}
\begin{footnotesize}
\noindent{\it Acknowledgement.}
We would like to thank especially T.\ Gannon for drawing our
attention to the non-symmetric $\SOs_1$ modular invariants
as well as K.-H.\ Rehren for pointing out the benefit of
proving Proposition \ref{extcft}, and we are indebted
to J.\ Fuchs, T.\ Gannon, and K.-H.\ Rehren for helpful
e-mail correspondences.
We gratefully acknowledge financial support of
the EU TMR Network in Non-Commutative Geometry.
\end{footnotesize}

%%%%%%%%%%%%%%%%%%%%%%%%%%%%%%%%%%%%%%%%%%%%%%%%%%%%%%%%%%%%%%%%%%%%%%%%%%%

%%%%%%%%%%%%%%%%%%%%%%%%% bibliography %%%%%%%%%%%%%%%%%%%%%%%%%%%%%%%%%%%%

\newcommand\bitem[2]{\bibitem{#1}{#2}}
%\marginpar{\scriptsize\tt{#1}}}

\def\aam              {Acta Appl.\ Math. }
\def\aip              {Ann.\ Inst.\ H.\ Poincar\'e (Phys.\ Th\'eor.) }
\def\cmp              {Com\-mun.\ Math.\ Phys. }
\def\duke             {Duke Math.\ J. }
\def\ijm              {Intern.\ J. Math. }
\def\jfa              {J.\ Funct.\ Anal. }
\def\jmp              {J.\ Math.\ Phys. }
\def\lmp              {Lett.\ Math.\ Phys. }
\def\rmp              {Rev.\ Math.\ Phys. }
\def\inv              {Invent.\ Math. }
\def\mpl              {Mod.\ Phys.\ Lett. }
\def\nup              {Nucl.\ Phys. }
\def\nupp             {Nucl.\ Phys.\ (Proc.\ Suppl.) }
\def\adma             {Adv.\ Math. }
\def\physa            {Physica {\bf A} }
\def\ijmp             {Int.\ J.\ Mod.\ Phys. }
\def\jp               {J.\ Phys. }
\def\fdp              {Fortschr.\ Phys. }
\def\plb              {Phys.\ Lett.\ {\bf B}}
\def\rims             {Publ.\ RIMS, Kyoto Univ. }

%%%%%%%%%%%%%%%%%%%%%%%%%%%%%%%%%%%%%%%%%%%%%%%%%%%%%%%%%%%%%%%%%%%%%%%%%%%

\begin{footnotesize}

\end{footnotesize}
\end{document}